\documentclass[review]{elsarticle}

\usepackage[utf8]{inputenc}
\usepackage{makeidx}         
\usepackage{graphicx} 
\usepackage{amssymb}
\usepackage{amsmath}
\usepackage{subcaption}
\usepackage{multicol}
\usepackage{algorithm, algorithmic}

\usepackage{todonotes}
\usepackage{layouts}
\usepackage[section]{placeins}

\usepackage{amsfonts}
\usepackage{color}

\newtheorem{exercise}{Exercise}
\newtheorem{exampl}{Example}
\newcommand{\comment}[1]{} 

\def\bq{\begin{quotation}}
\def\eq{\end{quotation}}

\definecolor{darkgreen}{rgb}{0.2,0.5,0.2}

\def\C{\Gamma}

\def\e{\varepsilon}

\def\F{\Phi}
\def\x{\chi}
\def\X{{\cal X}}
\def\y{\psi}
\def\Y{\Psi}

\def\O{\Omega}

\newcommand{\V}[1]{ \mathbf{#1} }    

\newcommand{\Vo}{\V{0}}

\newcommand{\Vb}{\V{b}}

\newcommand{\Ve}{\V{e}}
\newcommand{\Vf}{\V{f}}

\newcommand{\Vh}{\V{h}}

\newcommand{\Vn}{\V{n}}

\newcommand{\Vr}{\V{r}}

\newcommand{\Vt}{\V{t}}
\newcommand{\Vu}{\V{u}}
\newcommand{\Vv}{\V{v}}
\newcommand{\Vw}{\V{w}}

\newcommand{\Vy}{\V{y}}
\newcommand{\Vz}{\V{z}}

\newcommand{\Vgy}{\boldsymbol{\y}}

\newcommand{\M}[1]{ \mathbf{#1} }  

\newcommand{\MB}{\M{B}}
\newcommand{\MC}{\M{C}}

\newcommand{\MG}{\M{G}}
\newcommand{\MH}{\M{H}}
\newcommand{\MI}{\M{I}}

\newcommand{\MK}{\M{K}}

\newcommand{\MO}{\M{O}}

\newcommand{\MQ}{\M{Q}}
\newcommand{\MR}{\M{R}}
\newcommand{\MS}{\M{S}}
\newcommand{\MT}{\M{T}}
\newcommand{\MU}{\M{U}}
\newcommand{\MV}{\M{V}}
\newcommand{\MW}{\M{W}}

\newcommand{\MgC}{\M{\C}}
\newcommand{\MgF}{\M{\F}}

\newcommand{\MgY}{\M{\Y}}
\newcommand{\MgO}{\M{\O}}

\newcommand{\MgX}{\boldsymbol{\X}}

\newcommand{\Rmn}[2]{\mathbb{R}^{#1 \times #2}}

\def\eqv{\Leftrightarrow}

\def\2nm#1{\|#1\|_2}

\def\Ra#1{\mathrm{\bf range}(#1)}

\newcommand{\ars}[1]{\left[ \begin{array}{#1}}
\newcommand{\are}{\end{array} \right] }
\newcommand{\oars}[1]{\begin{array}{#1}}
\newcommand{\oare}{\end{array}}
\newcommand{\rars}[1]{\left( \begin{array}{#1}}
\newcommand{\rare}{\end{array} \right) }

\newcommand{\eqs}{\begin{eqnarray}}
\newcommand{\eqe}{\end{eqnarray}}
\newcommand{\eqsn}{\begin{eqnarray*}}
\newcommand{\eqen}{\end{eqnarray*}}

\def\defs{\begin{definition}}
\def\defe{\end{definition}}
\def\teos{\begin{theorem}}
\def\teoe{\end{theorem}}
\def\prfs{\begin{proof}}
\def\prfe{\end{proof}}
\def\exas{\begin{exampl}}
\def\exae{\end{exampl}}
\def\excs{\begin{exercise}}
\def\exce{\end{exercise}}
\def\cors{\begin{corollary}}
\def\core{\end{corollary}}

\newcommand{\ens}{\begin{enumerate}}
\newcommand{\ene}{\end{enumerate}}

\newcommand{\its}{\begin{itemize}}
\newcommand{\ite}{\end{itemize}}

\newcommand{\des}{\begin{description}}
\newcommand{\dee}{\end{description}}

\def\ul{\underline}
\def\wh{\widehat}

\def\wt{\widetilde}

\graphicspath{{figs/}}

\begin{document}
\begin{frontmatter}

\title{Krylov Subspace Recycling for Evolving Structures \tnoteref{mytitlenote}}

\author{M. Bolten}
\ead{bolten@math.uni-wuppertal.de}
\address{Faculty of Mathematics und Natural Sciences, University of Wuppertal, 42119 Wuppertal, Germany}

\author{E.~de Sturler\corref{corauthor}}
\ead{sturler@vt.edu}
\address{Department of Mathematics, Virginia Tech, Blacksburg, Virginia, 24061, USA}
\cortext[corauthor]{Corresponding Author}

\author{C. Hahn}
\ead{hahn@math.uni-wuppertal.de}
\address{Faculty of Mathematics und Natural Sciences, University of Wuppertal, 42119 Wuppertal, Germany}

\begin{abstract}
Krylov subspace recycling is a powerful tool when solving a long series of large, sparse linear systems that change only slowly over time. In PDE constrained shape optimization, these series appear naturally, as typically hundreds or thousands of optimization steps are needed with only small changes in the geometry. In this setting, however, applying Krylov subspace recycling can be a difficult task. As the geometry evolves, in general, so does the finite element mesh defined on or representing this geometry, including the numbers of nodes and elements and element connectivity. This is especially the case if re-meshing techniques are used. 
As a result, the number of algebraic degrees of freedom in the system changes, and in general the linear system matrices resulting from the finite element discretization change size from one optimization step to the next. Changes in the mesh connectivity also lead to structural changes in the matrices. In the case of remeshing, even if the geometry changes only a little, the corresponding mesh might differ substantially from the previous one. Obviously, this prevents us from any straightforward mapping of the approximate invariant subspace of the linear system matrix (the focus of recycling in this paper)
from one step to the next; similar problems arise for other selected subspaces.  
In this paper, we present an algorithm to map an approximate invariant subspace of the linear system matrix for the previous optimization step to an approximate invariant subspace of the linear system matrix for the current optimization step, for general meshes. This is achieved by exploiting the map from coefficient vectors to finite element functions on the mesh, combined with interpolation or approximation of functions on the finite element mesh. In addition, we develop a straightforward warm-start adaptation of the Krylov-Schur algorithm \cite{stewart2001} to improve the approximate invariant subspace at the start of a new optimization step if needed.
We demonstrate the effectiveness of our approach numerically with several proof of concept studies for a specific meshing technique. 
\end{abstract}

\begin{keyword}
\texttt{Moving Structures, Krylov Subspace Recycling, Recycling MINRES, Shape Optimization}
\MSC[2010] 00-01\sep  65-10
\end{keyword}

\end{frontmatter}

\section{Introduction}
\label{sec:intro}

In PDE constrained shape or topology optimization, usually hundreds or thousands of iterations are needed to reach a converged solution. Except for the first few steps, in many applications most of the changes in the shape per optimization step are small. In principle, this leads to sequences of matrices with only small changes in each optimization step. These correspond to small changes in the governing PDE, as well. Therefore, we expect that (preconditioned) Krylov subspace recycling may significantly speed up the computations, as exploited in \cite{wang2007}. In this paper, we focus on recycling approximate invariant subspaces.

However, in the shape optimization application presented here a significant obstacle to Krylov recycling arises as re-meshing is often necessary. This leads to matrices from one optimization step to the next that have different dimensions, which makes the linear systems algebraically incompatible regarding their sizes. This problem precludes a straightforward application of Krylov subspace recycling, since we cannot transfer the algebraic basis vectors spanning an approximate invariant subspace (or spanning any another useful subspace to recycle) from one linear system to the next. However, we will show that through the sequence of finite element spaces from which the algebraic linear systems are derived, we can map approximate algebraic eigenvectors from one linear system to the next by mapping (through interpolation or other approximations) the corresponding approximate eigenvectors, or basis vectors of an approximate invariant subspace, from one finite element space to the next. In a straightforward way, on AMR meshes (hence nested meshes), this idea was explored in \cite{WangPhD2007}; we consider more general mappings and corrections in this paper.
We will introduce the mapping for a generic mesh in Section \ref{sec:rminEG}.
This mapping through the finite element spaces may lead to additional approximation errors in the recycled subspace beyond the usual changes in compatible matrices. To improve the recycled approximate invariant subspace, when necessary, we also consider/propose a variation of the Krylov-Schur algorithm \cite{stewart2001} that has been adapted to allow a warm-start with a general subspace that is not a Krylov space for the new matrix. This Krylov-Schur variant is discussed in Section \ref{sec:rminEG-kryschur}. 
We choose an Arnoldi-based approach here, since RMINRES already includes a (cheap) Lanczos-based eigensolver. In addition, it allows the use of any preconditioner.
To show numerical examples, we give an example of a mapping for a specific meshing technique in Section \ref{sec:map}, the corresponding numerical examples are given in Section \ref{sec:num}.

\section{Foundations}
\label{sec:methods}
In shape optimization, the aim is to optimize a given two or three dimensional shape without changes in its inherent topology in order to maximize or minimize a given objective, for example efficiency or reliability. Depending on the problem considered, certain material laws have to be fulfilled, which are usually modeled by PDEs, such as the elasticity equations in structural mechanics or the heat equation. This leads, in principle, to a PDE constrained optimization problem, where the variable over which one optimizes is neither a function nor an element of a finite dimensional space, but a domain, i.e., in most cases a subset of $\mathbb{R}^2$ or $\mathbb{R}^3$. In this paper we restrict ourselves to domains in $\mathbb{R}^2$.
In a more formal way, let $\Omega$ be a bounded domain with Lipschitz boundary $\partial\Omega$ contained in a bounded set $\widetilde{\Omega}\subset\mathbb{R}^2$ such that each such $\Omega\subseteq\widetilde{\Omega}$ is admissible. This domain $\Omega$ is to be optimized subject to some functional $J(\Omega,u)$ that depends on $\Omega$ and the solution of a PDE $u\in H^1(\Omega)$, i.e. 
\begin{equation}
\begin{aligned}
&\!\min_{\Omega}        &\qquad& J(\Omega,u)\\
&\text{s.t.} &      & L u = f.
\end{aligned}
\end{equation}
We restrict ourselves to changes in the geometry that can be described by a mapping $F_t:\Omega\rightarrow \mathbb{R}^2$ that is a perturbation of identity, that is,
\begin{align*}
    F_t = \text{id}+t\mathcal{V},\quad t>0,
\end{align*}
where $\mathcal{V}\in (H^{1,\infty}(\Omega))^2$ is a velocity field.

\subsection{Krylov Subspace Recycling}
\label{sec:methods-recycl}
We briefly describe the main ideas behind Krylov subspace recycling and the recycling MINRES algorithm \cite{mottamello2010, wang2007}, a recycling version of the MINRES algorithm \cite{PaigeSaunders1975},  which is based on the earlier version from
\cite{KilmerdeSturler2006}. The purpose of recycling is to reuse a judiciously selected subspace from earlier linear solves to speed-up the convergence of subsequent linear solves. Here, we focus on approximate invariant subspaces associated with small (in absolute value) eigenvalues. In addition, depending whether the right hand sides change substantially or not, recycling may also give good initial guesses \cite{KilmerdeSturler2006}.

For alternative Krylov subspace recycling methods, for a range of applications, and HPC implementations, see \cite{jolivet2016block,Ahuja.Recyc-BiCGStab-MOD.2015,feng2013subspace,Gosselet.etal.reuse-Kryl-nonlin.2012,CGLV.FGCRODR.2012,Carlberg2015,Ahuja.Sturler.rGCROT-rBiCGStab-Hybrid.2015,al2018recycling}. A survey of Krylov subspace recycling methods is given in \cite{SoodStuKil2020}.  

In this section, we do not consider the complication that our matrices from one optimization step to the next may be incompatible. We will address this issue in section \ref{sec:rminEG}.

Consider at some optimization step the linear system $\MK \Vu = \Vf$ with $\MK \in \Rmn{n}{n}$, and let the space to be recycled be defined as $\Ra{\MW_k}$ with $\MW_k \in \Rmn{n}{k}$ computed in the previous optimization step. We compute $\wt{\MC}_k = \MK \MW_k$ and its thin QR decomposition $\wt{\MC}_k = \MC_k \MR_k$. Furthermore, assume an initial guess $\wt{\Vu}_0$ and corresponding residual $\wt{\Vr}_0 = \Vf - \MK \wt{\Vu}_0$. We first update the initial guess, adding the optimal correction (in minimum residual norm) from the space $\Ra{\MW_k}$.
We set $\Vu_0 = \wt{\Vu}_0 + \MW_k (\MR_k^{-1} \MC_k^T \wt{\Vr}_0)$, which gives the updated residual 
$\Vr_0 = \wt{\Vr}_0 - \MC_k \MC_k^T \wt{\Vr}_0 = (\MI - \MC_k \MC_k^T)\wt{\Vr}_0$.

We implement recycling Lanczos as follows. Let $\Vv_1 = \Vr_0 / \| \Vr_0 \|_2$.
We use the following augmented three-term recurrence
\begin{eqnarray}
\label{eq:Lanc_rec-1}
  \Vv_2 t_{2,1} & = & \MK\Vv_1 - \MC_k\Vb_1 - \Vv_1 t_{1,1}, \\
\label{eq:Lanc_rec-2}
  \Vv_{j+1}t_{j+1,j} & = & \MK\Vv_{j} - \MC_k\Vb_j - \Vv_{j}t_{j,j} - \Vv_{j-1} t_{j-1,j},
\end{eqnarray}
where $t_{j,j} = \Vv_j^T \MK\Vv_{j}$,
$\Vb_j = \MC_k^T \MK \Vv_j$,
$t_{j+1,j} = \| \MK\Vv_{j} - \MC_k\Vb_j - \Vv_{j}t_{j,j} - \Vv_{j-1} t_{j-1,j} \|_2$,
and $t_{j-1,j} = t_{j,j-1}$ was defined in the previous iteration,
leading to the augmented  Lanczos relation (with
$\MT$ for tridiagonal)
\eqs \label{eq:Lanc_rec-3}
  \MK \MV_j & = & \MC_k \MB_j + \MV_{j+1}\underline{\MT}_j 
  =
  \MC_k \MB_j + \MV_j \MT_j + \Vv_{j+1} \Ve_j^T t_{j+1,j} ,
\eqe
where $\MB_j = [\Vb_1 \;\; \Vb_2 \;\; \dots \;\; \Vb_j]$.
Note that the operator used in the Lanczos recurrence, 
$(\MI - \MC_k\MC_k^T)\MK$, is self-adjoint over the space $\Ra{\MC_k}^{\perp}$, which contains the 
space spanned by the Lanczos vectors, i.e.,
$(\MI - \MC_k\MC_k^T)\MK$ acts like a symmetric matrix on the vectors in the Krylov subspace  \cite{SoodStuKil2020}.
The approximate
solution in step $j$ of the recycling MINRES algorithm
is given by $\Vu_j = \Vu_0 + \MW_k \Vz + \MV_j \Vy$, 
where $\Vy$ and  $\Vz$ are determined by the minimum
residual condition. We have
\eqs
  \Vr_j = \Vf - \MK\Vu_j & = &
  \Vr_0 - \MC_k(\MR_k\Vz) - \MC_k\MB_j\Vy - 
  \MV_{j+1}\ul{\MT}_{j}\Vy \\
\label{eq:min_res-1}
  & = &
  \MV_{j+1}\left( \Ve_1 \|\Vr_0\|_2 - \ul{\MT}_{j}\Vy \right)
  - \MC_k \left( \MB_j\Vy + \MR_k\Vz \right)
\eqe
Due to the orthogonality $\MC_k^T \MV_{j+1} = \MO$, the minimization 
can be done in two separate steps. 
The first step, solving for $\Vy$, minimizes the residual over the space $\Ra{\MV_{j+1}}$, which is independent of the minimization over the space $\Ra{\MC_k}$. The minimization over $\Ra{\MC_k}$ simply requires that 
$\MB_j\Vy + \MR_k\Vz = \Vo$, and hence that
$\Vz = -\MR_k^{-1} \MB_j \Vy$. Note that the Lanczos recurrence is
the same as for standard MINRES, and we use a similar change of basis to develop a short term recurrence \cite{PaigeSaunders1975}. 
Using the thin QR decomposition of the 
tridiagonal matrix $\ul{\MT}_j = \ul{\MQ}_j^{(j+1)\times j} \MS_j^{j\times j}$, which is computed one column at a time, we recursively define 
$\wt{\MV}_j$ through $\wt{\MV}_j\MS_j = \MV_j$ (which leads to an additional three
term recurrence, as for MINRES), set $\wt{\Vy} =  \ul{\MQ}_j^T\Ve_1 \|\Vr_0\|_2$, 
and $\Vz = -\MR_k^{-1}\MB_j\MS_j^{-1} \wt{\Vy}$. 
We use the $\wt{\Vv}_j$ vectors to update the solution $\Vu_j$, so that the 
$O(n)$ vectors $\Vv_j$ and $\wt{\Vv_j}$ can be discarded. For
efficiency, the solution update $\MW_k \Vz$ can be postponed until after convergence. If the number of iterations is large, $\Vz$ can be updated recursively using an additional recurrence for the $O(k)$ vectors 
$\wt{\MB}_j \MS_j = \MB_j$. For details, see \cite{wang2007,mottamello2010}.

Algorithm~\ref{alg:RMINRES} outlines the Recycling MINRES algorithm that includes the recycle space into
the search space. For details on updating the
recycle space, see \cite{wang2007,mottamello2010}.
\begin{algorithm} 
\caption{Recycling MINRES}
\label{alg:RMINRES}
\begin{algorithmic}[1]
\STATE $\wt{\Vr}_0 = \Vf - \MK \wt{\Vu}_0$
\STATE $\Vu_0 = \wt{\Vu}_0$; \quad $\wh{\Vb} =  \MR_k^{-1}(\MC_k^T\wt{\Vr}_0)$; $\quad$
    $\Vr_0 = \wt{\Vr}_0-\MC_k\MC_k^T\wt{\Vr}_0$
\STATE $\Vv_1 = \Vr_0/\|\Vr_0\|_2$; $\quad$
    $\wt{\Vy} = \Ve_1 \|\Vr_0\|_2$
\FOR {$j=1, 2, \ldots$}
  \STATE $\wh{\Vv} = \MK \Vv_j$
  \STATE $\wh{\Vv} = \wh{\Vv} -\MC_k(\MC_k^T\wh{\Vv})$;
         $\quad$
         $\Vb_j = \MR_k^{-1}(\MC_k^T\wh{\Vv})$ \\
  \COMMENT{use modified Gram-Schmidt orthogonalization for updating $\wh{\Vv}$}
  \STATE $t_{j-1,j} =  t_{j,j-1}$; $\quad$
         $\wh{\Vv} = \wh{\Vv} - \Vv_{j-1}t_{j-1,j}$
  \STATE $t_{j,j} = \Vv_j^T \wh{\Vv}$; $\quad$
         $\wh{\Vv} = \wh{\Vv} - \Vv_j t_{j,j}$
  \STATE $t_{j+1,j} = \|\wh{\Vv}\|_2$; $\quad$
         $\Vv_{j+1} = \wh{\Vv}/t_{j+1,j}$
  \STATE $\MS_{:,j} = \MG_{j-1}\MG_{j-2}\ul{\MT}_{:,j}$ \\
  \COMMENT{apply the previous two Given's rotations to
    the new column of $\ul{\MT}_j$}
  \STATE Compute Given's rotation $\MG_j$ such that
         $\MS_{:,j} = \MG_j\MS_{:,j}$ has
         $s_{j+1,j} = 0$.
  \COMMENT{see MINRES \cite[p. 41--44]{GreenbaumBk1997}}
  \STATE $\wt{\Vy} = \MG_j\wt{\Vy}$
  \STATE $\wt{\Vv}_j = s_{j,j}^{-1}(\Vv_j - \wt{\Vv}_{j-1}s_{j-1,j} - \wt{\Vv}_{j-2}s_{j-2,j})$
  \STATE $\wt{\Vb}_j =  s_{j,j}^{-1}(\Vb_j - \wt{\Vb}_{j-1}s_{j-1,j} - \wt{\Vb}_{j-2}s_{j-2,j})$
  \STATE $\Vu_j = \Vu_{j-1} + \wt{\Vv}_j\wt{y}_j$
   \qquad \COMMENT{$\wt{y}_j$ is the $j$th entry of vector $\wt{\Vy}$}
  \STATE $\wh{\Vb} = \wh{\Vb} - \wt{\Vb}_j\wt{y}_j$ \\
  \STATE {\bf if} $| \wt{y}_{j+1} | \leq \e \|\Vr_0\|_2$ {\bf then }{\tt break}
\ENDFOR
\STATE $\Vu = \Vu_j + \MW_k \wh{\Vb}$
\end{algorithmic}
\end{algorithm}

While algorithm \ref{alg:RMINRES} allows any subspace to be recycled, we focus here on approximate invariant subspaces, in particular, those corresponding to small eigenvalues. First, this leads to substantially improved rates of convergence as the condition number is significantly improved. Second, the small eigenvalues correspond to smooth modes that can be transferred effectively from one iteration of the shape optimization to the next.

\section{Recycling MINRES for evolving geometries}
\label{sec:rminEG}

In shape optimization, the changes in geometry in each optimization step and thus the underlying mesh prevent a straightforward application of the described Krylov subspace recycling. Depending on the meshing technique, a mapping of the matrix representing the subspace in one optimization step to the next might be necessary. 

\subsection{Mapping between successive meshes}
\begin{figure}
    \centering
    \includegraphics[width=0.95\textwidth]{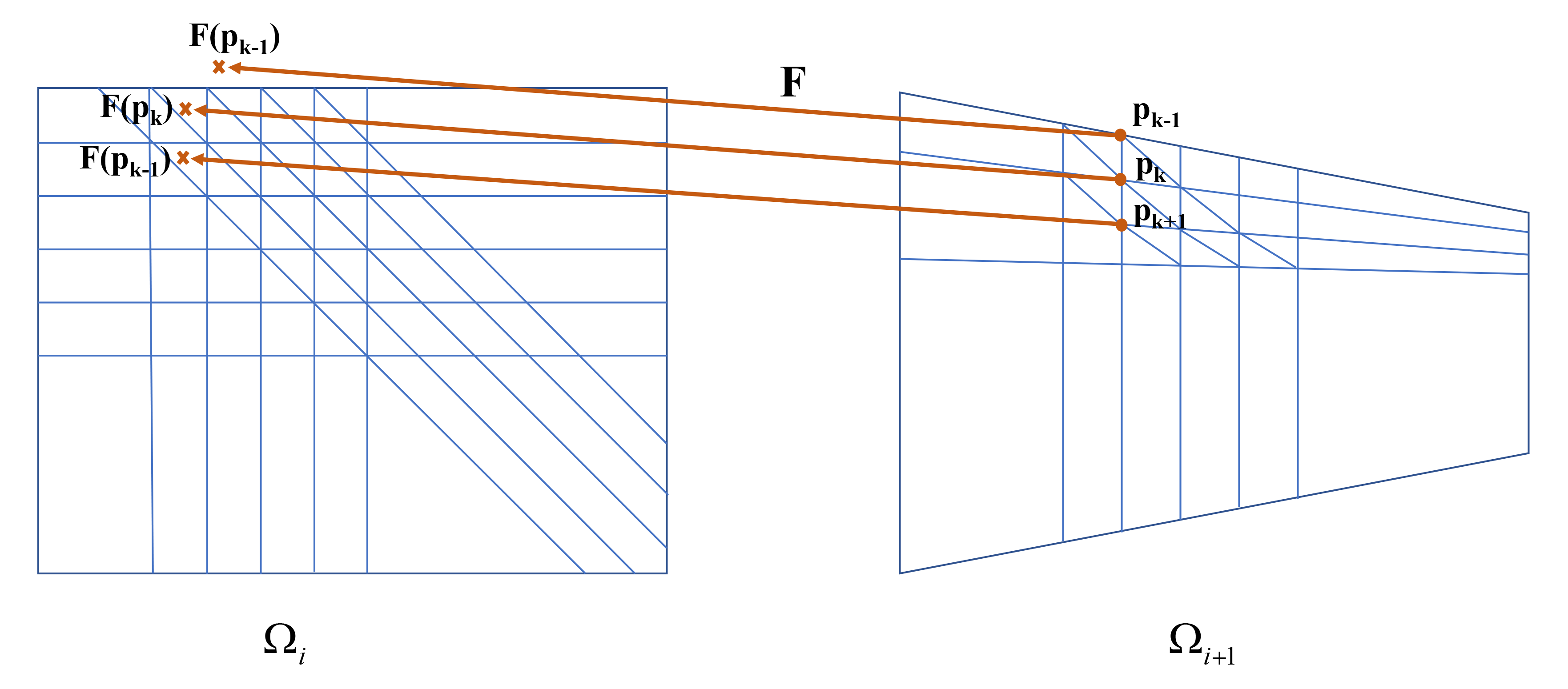}
    \caption{Mapping new mesh nodes to the domain of the previous optimization step.}
    \label{fig:dif-shapes}
\end{figure}

Let $\Omega_i$ and $\Omega_{i+1}$ be two domains representing two immediately consecutive shapes stemming from an iterative shape optimization procedure, see Figure \ref{fig:dif-shapes}. These domains are discretized by finite element meshes $T_i$ and $T_{i+1}$, respectively. In general, the meshes feature different connectivities and different numbers of nodes, $N_i$ and $N_{i+1}$, respectively.
Additionally, from optimization step $i$, an approximate invariant subspace is given by $\Ra{\MW_k^{(i)}}$ with $\MW_k^{(i)} \in \Rmn{N_i}{k}$, that is supposed to be recycled in optimization step $i+1$. The system in this step however is of dimension $N_{i+1}$. Therefore, we are in need of a function that maps the ${N_i}\times{k}$ - dimensional matrix representing the approximate invariant subspace of the system in optimization step $i$ to a $N_{i+1}\times{k}$ - dimensional matrix, that hopefully still represents a good approximate invariant subspace of the linear system in optimization step $i+1$.

The mapping we propose exploits the fact that the algebraic systems we consider are closely linked to the continuous finite element spaces $\mathcal{V}_h(\Omega_i, T_i)$ and $\mathcal{V}_h(\Omega_{i+1}, T_{i+1})$, determined by $T_i$ and $T_{i+1}$, respectively. As a reminder, we give a recap to the Galerkin method: In general, to solve a PDE, it is considered in its weak formulation defined on a Sobolev space $H^m(\Omega)$ or $H_0^m(\Omega)$. The solution is approximated in a finite dimensional subspace of this Sobolov space, the finite element space $\mathcal{V}_h(\Omega)$. Consider for example the elliptic zero-boundary value problem 
\begin{align}
\begin{split}
    -\sum\limits_{i,k=1}^n\partial_i(a_{ik}\partial_k u) + a_0 u &= f \quad \text{in } \Omega\\
    u &= 0 \quad \text{on } \partial\Omega,
    \end{split}
    \label{eq:elip}
\end{align}
with $u\in H_0^1(\Omega)$. The function $u$ is a solution to (\ref{eq:elip}) if
\begin{align}
    a(u,v) = \langle f,v\rangle, \quad \text{for all } v \in H_0^1(\Omega).
    \label{eq:galerkin}
\end{align} 
On the finite dimensional subspace $V_h(\Omega)$ we can find a basis of nodal basis functions $\Phi_1,\dots,\Phi_N$ of the subspace such that (\ref{eq:galerkin}) is equivalent to 
\begin{align}
    a(u_h,\Phi_j) = \langle f,\Phi_j\rangle, \quad j=1,\dots,N,
    \label{eq:galerkin_finite1}
\end{align} 
with $u_h\in V_h(\Omega)$. 
Therefore, $u_h$ can be formulated in terms of the basis functions, $u_h =\sum_{j=1}^Nu_j\Phi_j$. This approach leads to the system of equations 
\begin{align}
    \sum\limits_{k=1}^N a(\Phi_k,\Phi_j)u_j = \langle f,\Phi_j\rangle, \quad j=1,\dots,N,
    \label{eq:galerkin_finite2}
\end{align} 
or in matrix notation $\MK\Vu=\Vb$.
Now, consider the approximate invariant subspace $\Ra{\MW_k^{(i)}}$ with $\MW_k^{(i)} \in \Rmn{N_i}{k}$. Instead of considering the coefficients of the matrix in the algebraic sense, we can see the columns of the matrix as vectors containing the coefficients of continuous functions defined on the finite element space $\mathcal{V}_h(\Omega_i, T_i)$. These functions are defined as 
\begin{align}
    w_m^{(i)}(x):= \sum\limits_{j=1}^{N_i}(\MW_k^{(i)})_{(j,m)}\Phi_j^{(i)}(x), \quad m=1,\dots,k,\quad x\in\Omega_i,
    \label{eq:map-matrix}
\end{align}
$\Phi_j^{(i)}(x)$ being the nodal basis functions of the finite element space $\mathcal{V}_h(\Omega_i, T_i)$. 
If we are able to find a transformation $F(x^{(i+1)}):\Omega_{i+1}\rightarrow\Omega_i$, that maps each node $p_l$ in $T_{i+1}$, $l=1,\dots,N_{i+1}$ to a corresponding point $x_l$ in $\Omega_i$, we can build the matrix representing the mapped approximate invariant subspace, $\widetilde{\MW}_k^{(i)} \in \Rmn{N_{i+1}}{k}$, in the following way:
\begin{align}
    (\widetilde{\MW}_k^{(i)})_{(l,m)} := w_m^{(i)}(F(p_l)).
\end{align}
In summary, this method consists of two steps: First, we define a map from the nodes of the new mesh to points in the old domain (see Fig. \ref{fig:dif-shapes}); and second, we interpolate the values of the matrix $\MW_k^{(i)}$ via the functions defined in (\ref{eq:map-matrix}) in the finite element space and evaluate these functions at the points in $\Omega_i$ corresponding to the nodes in of $T_{i+1}$.

Although the second step is straightforward, the first step can be challenging, depending on the meshing technique that is used. Therefore, we briefly discuss general concepts to realize such a mapping for different meshing techniques.
If the update of the shape in the optimization procedure is performed via \textbf{mesh morphing} \cite{staten2011}, each node in $T_{i+1}$ can be mapped uniquely to the corresponding node in $T_i$, hence the approximate invariant subspace can be recycled without being transformed. Techniques that work with a reference frame, like \textbf{Arbitrary Lagrangian Eulerian (ALE)} \cite{hirt1974}, come with an inherent mapping from the reference frame to the domain. Through this reference frame, we can map each node in $T_{i+1}$ to a node in $T_i$ and thus recycle the approximate invariant subspace as in the mesh morphing case. For meshing techniques that do not have such an inherent mapping it gets more difficult.

If we assume, for example, a general re-meshing scheme without a prescribed number of nodes or restrictions on the connectivity, then the optimization procedure does not provide any information about the relation between the two domains $\Omega_i$ and $\Omega_{i+1}$; see Figure  \ref{fig:dif-shapes}. Additionally, the two meshes generally may differ in the number of nodes and thus basis functions. Therefore, a more sophisticated  mapping of the approximate invariant subspace of the linear system matrix derived from mesh $T_i$ to 
an approximate invariant subspace for the system matrix derived from mesh $T_{i+1}$ becomes inevitable. 
If no other information is available, the simplest map from $\Omega_{i+1}$ to $\Omega_i$ is $F(x) := x$. This choice, however, does not guarantee that  $F(x)\in\Omega_i$. Therefore, we suggest to choose $F(x)$ as the minimizer, in a suitable norm, of the distance between the given point $x$ and points in $\Omega_i$, i.e., 
$F(x):= \arg \min_{\tilde{x}_i\in\Omega_i}
||\tilde{x}_i-x||$, which implies the identity for $x \in \Omega_i$. Another approach  
would be to extrapolate the $w_m^{(i)}(\cdot)$ in a simple way. We discuss this idea in more detail in section \ref{sec:map}. 

A special case of re-meshing is mesh refinement (and derefinement), especially \textbf{adaptive mesh refinement} (AMR) \cite{berger1984}. Although in this case the two systems will certainly differ in dimension, 
the new nodes will definitely lie inside $\Omega_i$. Additionally, for the new nodes the global evaluation in (\ref{eq:map-matrix}) reduces to a local evaluation on the respective refined element; a mapping of approximate invariant subspaces between AMR meshes is described in  \cite{WangPhD2007}.

The discussed approach results in two obvious challenges: The first is already addressed at the end of section \ref{sec:methods-recycl}. To reduce the condition number of a system, it sometimes can be desirable to remove the absolute largest eigenvalues rather than the smallest ones or both the smallest and the largest. Yet, the large eigenvalues often correspond to high frequency modes which are clearly very difficult to approximate by the presented approach. We therefore restrict ourselves to remove small eigenvalues, which generally is a good choice. The second challenge is that small changes in the geometry do not necessarily result in small changes in the eigenvectors of the system matrix. In such cases, the approximate invariant subspace of the system matrix for the new geometry and mesh, computed by the mapping described above, is not sufficiently accurate that it leads to a fast convergence rate. This issue is addressed in section \ref{sec:rminEG-kryschur}.

\subsection{A Warm-start Krylov-Schur Algorithm}
\label{sec:rminEG-kryschur}

While our approach to map an approximate invariant subspace from optimization step to the next, discussed in the previous subsection, usually provides good approximate invariant subspaces, the error from the approximation procedure is sometimes large enough that recycling the approximate invariant subspace is not effective. In such a case, we can use a few iterations of an eigensolver to improve the approximation.

The Krylov-Schur algorithm \cite{stewart2001} is a (more stable) variant of the Implicitly Restarted Arnoldi (IRA) algorithm \cite{Sore_IRA1_1992}; it builds an alternative Krylov recurrence
that takes the form
\eqs \label{eq:KrylSchur_1}
  \MK \MU_m & = & \MU_m \MS_m + \Vu_{m+1}\Vh_{m+1}^T = 
  \MU_{m+1} \ul{\MS}_m ,
\eqe
where $\MS_m$ is upper triangular and $\Vh_{m+1}$ is a dense vector (in principle), giving an approximate block Schur decomposition (apart from the residual term). The key idea of the Krylov-Schur algorithm is that after ordering the approximate eigenvalues on the diagonal in the desired order, the Krylov-Schur decomposition can be truncated at any point. Writing (\ref{eq:KrylSchur_1}) as the block decomposition 
\eqs\label{eq:KrylSchur_2}
  \MK \ars{cc} \MU^{(1)} & \MU^{(2)} \are  & = & 
    \ars{cc} \MU^{(1)} & \MU^{(2)} \are 
  \ars{cc} \MS_{11} & \MS_{12} \\ \MO & \MS_{22} \are \\
    && + \Vu_{m+1} 
    \ars{cc} \left( \Vh^{(1)}_{m+1} \right)^T & \left( \Vh^{(2)}_{m+1} \right)^T \are ,
\eqe
where the superscript $(1)$ indicates the desired eigenvalues and Schur vectors,
we can form the new, truncated, Krylov-Schur 
decomposition 
\eqs\label{eq:KrylSchur_3}
  \MK \MU^{(1)} & = & 
    \MU^{(1)} \MS_{11} + \, \Vu_{m+1} \left( \Vh^{(1)}_{m+1} \right)^{T} .
\eqe
Subsequently, we can expand the Krylov decomposition using additional iterations in the same way as for standard Arnoldi. After this, we can again put the recurrence in the form of (\ref{eq:KrylSchur_1}), truncate, and continue until the method converges.

Next, we describe a variant of Stewart's Krylov-Schur algorithm to update our approximate invariant subspace. However, starting with a set  of vectors that do not form a Krylov sequence or span a Krylov space leads to several complications in the algorithm. 

The first complication is that we need to compute the best approximate Krylov
decomposition \cite{stewart2002}. 
Let $\MK$ be the matrix in the current optimization step, and let $\wt{\MW} \in \Rmn{N_{i}}{k}$ be the matrix whose columns span an approximate invariant subspace. Following \cite{stewart2002}, we compute the approximate Krylov decomposition (matrix) $[ \MU_{k-1} \; \Vu_{k} ]$ with 
$\Ra{[\MU_{k-1} \; \Vu_{k} ]} = \Ra{\wt{\MW}}$ with the smallest backward error for $\MK$. 
We start with the thin QR-decomposition $\wt{\MW} = \wt{\MU}_{k} \wt{\MS}_{k}$, set 
$\wt{\MH}_{k} = \wt{\MU}_{k}^T\MK \wt{\MU}_{k}$, and residual $\wt{\MR} = \MK\wt{\MU}_{k} - \wt{\MU}_{k} 
\wt{\MH}_{k}$. Next, we compute the 
matrix $\MgY$ in the thin singular value decomposition (SVD) 
$\wt{\MR} = \MgF \MgO \MgY^T$, and set 
$\MV_{k} = [\Vgy_2 \; \Vgy_3 \; \ldots \; \Vgy_{k}\; \Vgy_1]$ (hence $\MV_{k-1} = [\Vgy_2 \; \Vgy_3 \; \ldots \; \Vgy_{k}]$). This gives
\eqs
\label{eq:KrylDecomp_1}
\MK (\wt{\MU}_{k} \MV_{k-1}) & = &
  (\wt{\MU}_{k}\MV_{k}) (\MV^T_{k}\wt{\MH}_{k}\MV_{k-1}) + 
  \wt{\MR}\MV_{k-1}  \eqv \\
\label{eq:KrylDecomp_2}
\MK \MU_{k-1} & = &
  [\MU_{k-1} \;\; \Vu_{k}]\ul{\MH}_{k-1} +
  \MR , 
\eqe
where $\MU_{k-1} = \wt{\MU}_{k} \MV_{k-1}$, 
$\Vu_{k} = \wt{\MU}_{k} \Vv_{k} = \wt{\MU}_{k} \Vgy_1$, 
$\ul{\MH}_{k-1} = \MV^T_{k}\wt{\MH}_{k}\MV_{k-1}$, 
and $\MR = \wt{\MR}\MV_{k-1}$.
Following \cite{stewart2002}, 
(\ref{eq:KrylDecomp_2}) gives an approximate Krylov decomposition for $\MK$ and $\Ra{\wt{\MW}}$
with a minimal backward error (i.e., with minimal norm residual). 

Next, we start the Warm-start Krylov-Schur algorithm, extending the search space using $m-k+1$ Arnoldi steps. This gives,
\eqs \label{eq:KrylDecomp_Ext_1}
\MK \MU_{m} & = &
    \MU_{m+1} 
  \ars{cc}
    \ul{\MH}_{k-1} & \MH_{12} \\  
    \MO & \MH_{22} 
  \are
  + 
  \ars{cc} \MR &\vline \;\; \MO \are ,
\eqe
where 
$\MU_{m+1} = [\MU_{k} \;|\; \MU_c]
=[\Vu_1 \; \Vu_2 \; \ldots \; \Vu_{k} \;|\; \Vu_{k+1} \; \ldots \; \Vu_{m+1}]$ with 
$\MU_c = [\Vu_{k+1} \; \ldots \; \Vu_{m+1}]$.
We put this recurrence in quasi Krylov-Schur form to perform a truncation following 
(\ref{eq:KrylSchur_1})--(\ref{eq:KrylSchur_3}). However, since we did not start with a Krylov space, we need to account for the residual matrix in computing the Rayleigh quotient. 
Since 
$\Ra{\wt{\MU}_{k}} = \Ra{\MU_{k}}$, we have
$\MU_{k}^T \MR = \MO$ and 
\eqs \label{eq:CorrRaylQuo}
  \ul{\MH}_m & = & 
  \MU_{m+1}^T \MK \MU_{m} = 
  \ars{cc}
    \ul{\MH}_{k-1} & \MH_{21} \\
    \MO & \MH_{22} 
  \are + 
  \MU_{m+1}^T\ars{c} \MR \;\vline\; \MO \are \\
  & = & 
  \ars{cc}
    \ul{\MH}_{k-1}  &  \MH_{21} \\
    \MU_c^T \MR     &  \MH_{22} 
  \are ,
\eqe
where $\MU_c^T\MR$ can
be computed as 
$(\MU_c^T\wt{\MR})\MV_{k-1}$  
in the first Arnoldi iteration. 
This gives the new residual
\eqs \label{eq:New_Residual}
  \MR_\mathrm{new} & = &
  \MK \MU_m - \MU_{m+1} \ul{\MH}_m ,
\eqe
with 
$\MU_{m+1}^T \MR_\mathrm{new} = \MO$.
We note that only the leading $k-1$ columns of $\MR_\mathrm{new}$
can be nonzero. Moreover,
$\| \MR_\mathrm{new} \|_F \leq \|\MR\|_F \leq \| \wt{\MR} \|_F$, since we have 
from (\ref{eq:KrylDecomp_Ext_1}) and 
(\ref{eq:CorrRaylQuo}),
\eqs 
\nonumber
  \MR_\mathrm{new} & = & 
  \MK \MU_m - \MU_{m+1}\ul{\MH}_m \\
  & = &
  \ars{c} \MR \;|\: \MO \are -
  \MU_c\MU_c^T \ars{c} \MR \;|\: \MO \are .
\eqe

To put the recurrence in Krylov-Schur form, 
we take the real Schur decomposition of the leading 
$m \times m$ part of  $\ul{\MH}_m$, 
$\MgX \MgC \MgX^T = \MH_m$, with 
the desired eigenvalues for the truncation in $\MgC_1$, the leading
$(k-1) \times (k-1)$ (or if necessary $k \times k$) part of $\MgC$, 
and make a change
of basis using $\MgX$,
as follows.
\eqs
\label{eq:KrylDecomp_11}
\MK \MU_{m} & = &
  [\MU_{m} \;\; \Vu_{m+1}]
  \ars{c} \MgX \MgC \MgX^T \\ \Vh^T 
  \are  + \MR_\mathrm{new} \;\;\eqv \\
\label{eq:KrylDecomp_12}
\MK ( \MU_{m} \MgX) & = &
  [ (\MU_{m}\MgX) \;\; \Vu_{m+1} ]
  \ars{c} \MgC  \\ \Vh^T\MgX 
  \are  + \MR_\mathrm{new} \MgX.
\eqe
Next, we split the Schur decomposition into two parts, with the first part corresponding to the desired eigenvalues. 
\eqs
\nonumber
\MK [ \MU_{m}\MgX^{(1)} \;\; \MU_{m}\MgX^{(2)} ] & = &
  [\MU_{m}\MgX^{(1)} \;\; \MU_{m}\MgX^{(2)} \;\; \Vu_{m+1} ]
  \ars{cc} \MgC_{11}  & \MgC_{12} \\
           \MO        & \MgC_{22} \\ 
           \Vh^T\MgX^{(1)} & \Vh^T\MgX^{(2)} 
  \are \\  
  && + \; [\MR_\mathrm{new}\MgX^{(1)} \;\; \MR_\mathrm{new}\MgX^{(2)}] 
  \qquad \eqv 
\\
\nonumber
\MK [\MU^{(1)} \;\; \MU^{(2)} ] & = &
  [\MU^{(1)} \;\; \MU^{(2)} \;\; \Vu_{m+1} ]
  \ars{cc} \MgC_{11}  & \MgC_{12} \\
           \MO        & \MgC_{22} \\ 
           \Vh_1^T & \Vh_2^T 
  \are \\
  && + \; [\MR_\mathrm{new}^{(1)} \;\; \MR_\mathrm{new}^{(2)}] .
\label{eq:KrylDecomp_13}
\eqe
Finally, we truncate the Krylov-Schur 
decomposition with residual term, giving
\eqs \label{eq:KrylDecomp_14} 
\MK \MU^{(1)} & = &
  [\MU^{(1)} \;\; \Vu_{m+1} ]
  \ars{c} \MgC_1 \\ \Vh_1^T \are
  + \MR_\mathrm{new}^{(1)} .
\eqe

With $\MU_{k} = [\MU^{(1)} \;\; \Vu_{m+1} ]$
and $\MR = \MR_\mathrm{new}^{(1)}$ 
(\ref{eq:KrylDecomp_2}) is again satisfied, and we can repeat the 
Arnoldi-based extension and truncation, another Warm-Start Krylov-Schur cycle, until a desired accuracy is reached.

\begin{algorithm}
\caption{Warm-start Krylov-Schur}
\begin{algorithmic}[1]
\STATE Input: matrix $\MK$,
  basis $\MW_{k}  = [\Vw_1, \ldots, \Vw_{k}]$,
  nr. of cycles $J$, cycle length $m$, recycle space dimension $k$ \\
\COMMENT{{\em Initialize: Compute Krylov decomposition with min backward error}}
\STATE 
  $[\MU_{k},\MS_{k}] = 
  \mathrm{thinQR}(\MW_{k})$; $\:$
  $\MW_{k} = \MK \MU_{k}$; $\:$
  $\MH = \MU_{k}^T\MW_{k}$;
\STATE 
  $\MR = \MW_{k} - \MU_{k}\MH$;
  $[\MgF,\MgO,\MgY] = \mathrm{thinSVD}(\MR)$; 
  $\MV_{k} = [\Vgy_2 \, \ldots \, \Vgy_{k} \, \Vgy_1]$;
\STATE 
  $\MH = \MV_{k}^T(\MH \MV_{k-1})$; \,\,
  $\MU_{k} = \MU_{k}\MV_{k}$; \,\,
  $\MR = \MR \MV_{k-1}$; \\
\COMMENT{{\em End Initialize}} \\
\COMMENT{{\em $J$ restarts of Warm-start Krylov-Schur}}
\FOR {$j=1,\ldots, J$}
  \STATE{\bf if}{ space converged }
  {\bf then}{ break}
  {\bf end if} \\
  \COMMENT{{\em $m-k+1$ Arnoldi steps}}
  \FOR{$i = k, \ldots, m$} 
    \STATE 
    $\Vt = \MK \Vu_i$;
    \STATE 
    $\MH_{:,i} = \MU_i^T\Vt$; $\:$
    $\Vt = \Vt - \MU_i\MH_{:,i}$; \,\,
    \COMMENT{{\em accurate orthog. needed}}
    \STATE 
    $\MH_{i+1,i} = \|\Vt\|_2$; $\:$
    $\Vu_{i+i} = \Vt / \MH_{i+1,i}$;
  \ENDFOR \\
  \COMMENT{{\em Update Rayleigh quotient and residual}}
  \STATE 
    $\MH_{k+1:m+1,1:k-1} = \MU_{1:n,k+1:m+1}^T \MR$; 
  \STATE
    $\MR = \MR - \MU_{1:n,k+1:m+1} \MH_{k+1:m+1,1:k-1}$; \\
  \COMMENT{{\em compute Schur form and truncate}}
  \STATE 
    $[\MgX,\MgC] = \mathrm{realSchur}(\MH_{1:m,1:m})$; \,
    $\MU_{k-1} = \MU_m \MgX_{1:m,1:k-1}$; \,
    $\Vu_{k} = \Vu_{m+1}$;
  \STATE
    $\MH_{1:k-1,1:k-1} = 
      \MgC_{1:k-1,1:k-1}$; \, 
    $\MH_{k,1:k-1} = \MH_{m+1,1:m} \MgX_{1:m,1:k-1}$;
  \STATE $\MR = \MR \MgX_{1:k-1,1:k-1}$
\ENDFOR
\end{algorithmic}
\end{algorithm}

Some further efficiency improvements are possible. For
example, we can skip computing 
$\MU_{k} = \wt{\MU}_{k}\MV_{k}$ after the initialization.
In the first iteration of the Warm-start Krylov-Schur algorithm,
orthogonalization can be done with $\wt{\MU}_{k}$ and 
the corresponding columns of $\MH$ multiplied by $\MV_{k}^T$, since 
$\MU_{k}\MU_{k}^T \Vt = \wt{\MU}_{k}\wt{\MU}_{k}^T \Vt$, and $\MU_{k}^T \Vt = \MV_{k}^T \wt{\MU}_{k}^T \Vt$.

For the purposes of effective recycling in the linear solver, we generally need only a modest
improvement of the accuracy of the desired recycle space as an approximate invariant subspace
\cite{ParkStur2006}. Hence, we suggest a fixed few cycles of the Warm-Start Krylov-Schur algorithm.
In the numerical example in section \ref{sec:num-ks}, we see an improvement in the convergence already after two cycles of the Warm-Start Krylov-Schur algorithm.

\subsection{Cost of the Warm-Start 
Krylov-Schur Algorithm}
We split the 
Warm-Start Kylov-Schur algorithm into two steps:
(1) The computation of the approximate 
Krylov decomposition with minimum backward error, and (2) a few cycles of the Warm-Start Krylov-Schur algorithm. 
We need to compute the (thin) QR-decomposition of  $\wt{\MW}_i$ at (approximately) $2Nk^2$ flops, the Rayleigh quotient $\wt{\MH}_{k}$ and the 
residual $\wt{\MR}$ at the cost of $k$ sparse preconditioned matrix-vector products plus the 
cost of the orthogonalizations, (approximately) $4Nk^2$ flops, and the 
$\MgY$ factor of the singular value decomposition at (approximately) $2Nk^2$ flops.
The 
initial matrices $\MR$ and 
$\MU_{k}$ need not be explicitly computed. The total cost of this step, counting only terms including $Nk^2$ ($k^2$ being the leading power of $k$), 
is approximately $8Nk^2$ flops.

The first part of the Warm-start Krylov-Schur cycle
consists of $m-k+1$ Arnoldi iterations. This requires 
$m-k+1$ preconditioned matrix-vector products and
$(1/2)(m-k+1)(m+k)$ orthogonalizations; the latter takes
$2(m-k+1)(m+k)N$ flops. Updating the Rayleigh quotient takes an additional $2(m-k+1)(k-1)N$ for the inner products to compute $\MU_c^T \MR$. 
The Schur decomposition of the Rayleigh quotient does not require any $O(Nk^2)$
computations. After determining $\MgX^{(1)}$, if another Krylov Schur cycle is needed, computing the matrices costs approximately $2N(m+k-1)(k-1)$ flops. Note that if $m-k$ is relatively small there is a cheaper way to compute the basis to restart; see \cite{Stu99}. Summarizing, the total cost of one cycle of the Warm-start Krylov-Schur is 
$m-k+1$ preconditioned matrix-vector products plus $2N(m^2 - k^2+2mk)$ flops. A cycle of the Warm-start Krylov Schur algorithm is only slightly more expensive 
than a cycle of the standard Krylov Schur algorithm, but with a good starting space the Warm-Start version tends to converge much faster.

\section{Mapping on structured meshes}
\label{sec:map}

To solve the optimization problem given in section \ref{sec:methods}, it is first discretized via the finite element method. Here, the meshing is realized by an approach similar to Composite Finite Elements, first developed in \cite{hackbusch1996, hackbusch1997a, hackbusch1997b}.

\subsection{Structured meshes for evolving geometries}
\label{sec:methods-cfe}

The structured meshing approach that is used in this paper combines some of the best features of re-meshing and mesh morphing techniques. For demonstration purposes, the technique is described in two dimensions but easily extends to three dimensions. For more details see \cite{bolten2020}.

\begin{figure}[t]
	\centering
	\begin{subfigure}[c]{0.32\textwidth}
		\includegraphics[width=0.75\linewidth]{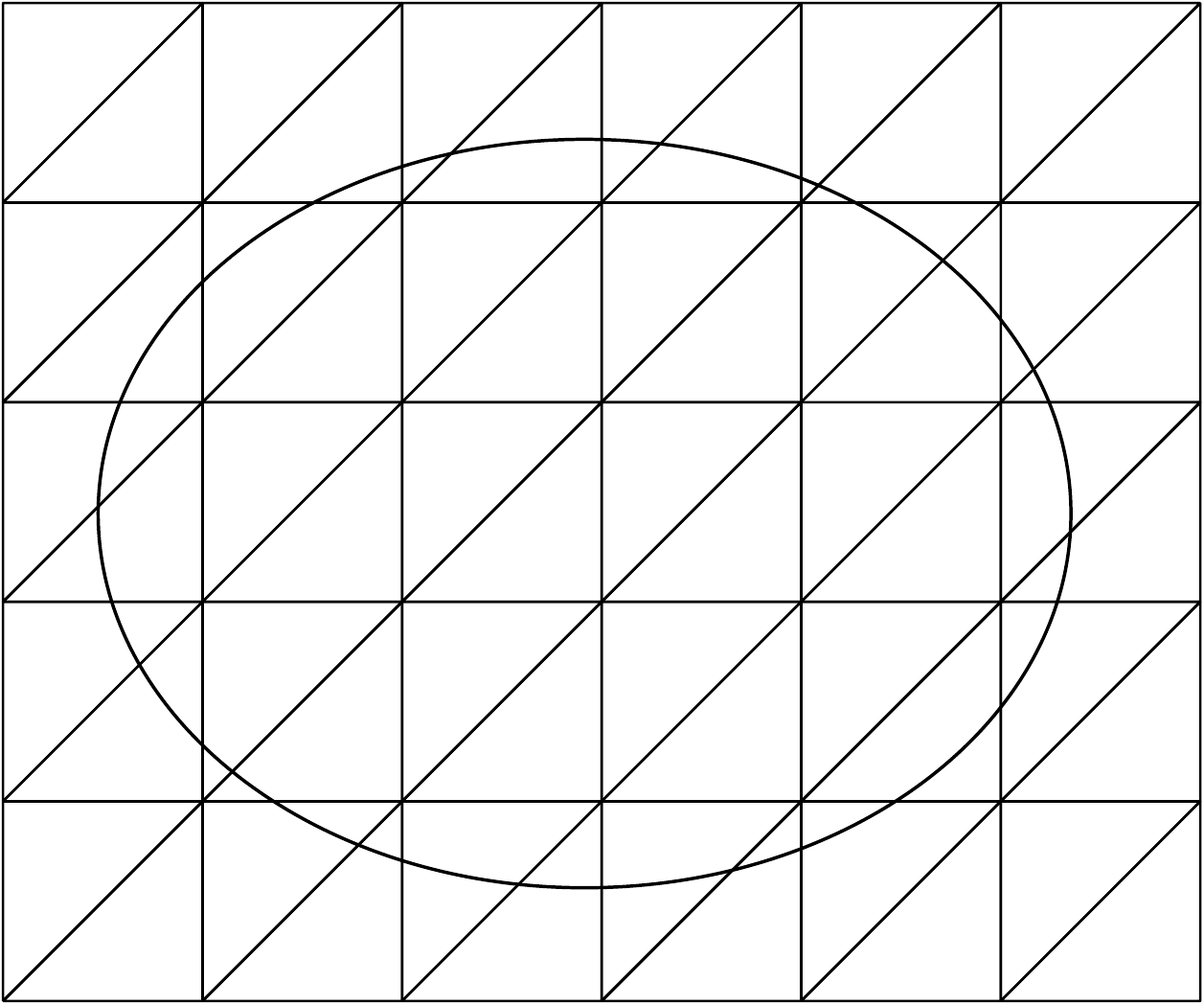}
		\caption{ }
		\label{fig:cfe-a}
	\end{subfigure}
	\begin{subfigure}[c]{0.32\textwidth}
		\includegraphics[width=0.75\linewidth]{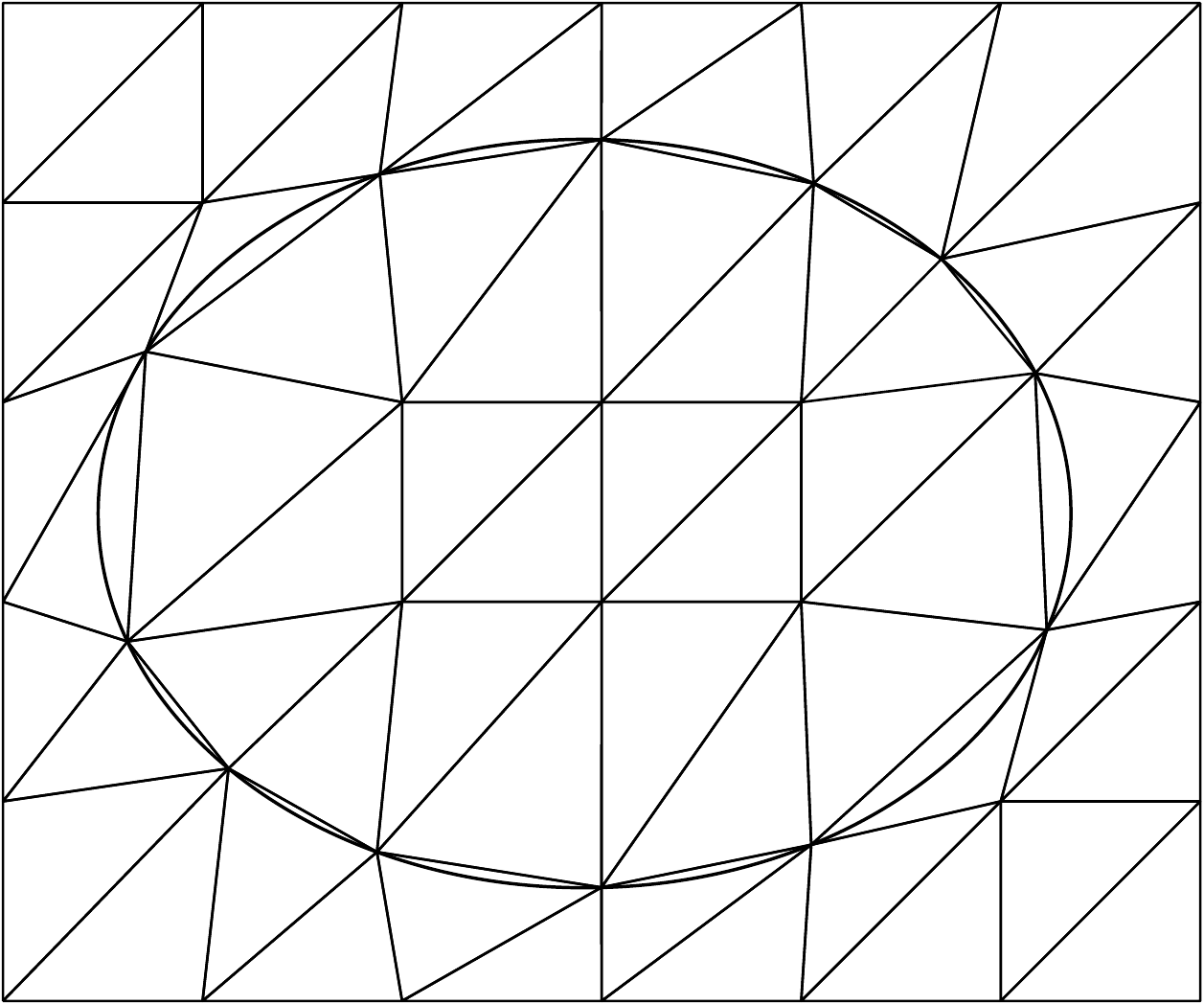}
		\caption{ }
		\label{fig:cfe-b}
	\end{subfigure}
	\begin{subfigure}[c]{0.32\textwidth}
		\includegraphics[width=0.75\linewidth]{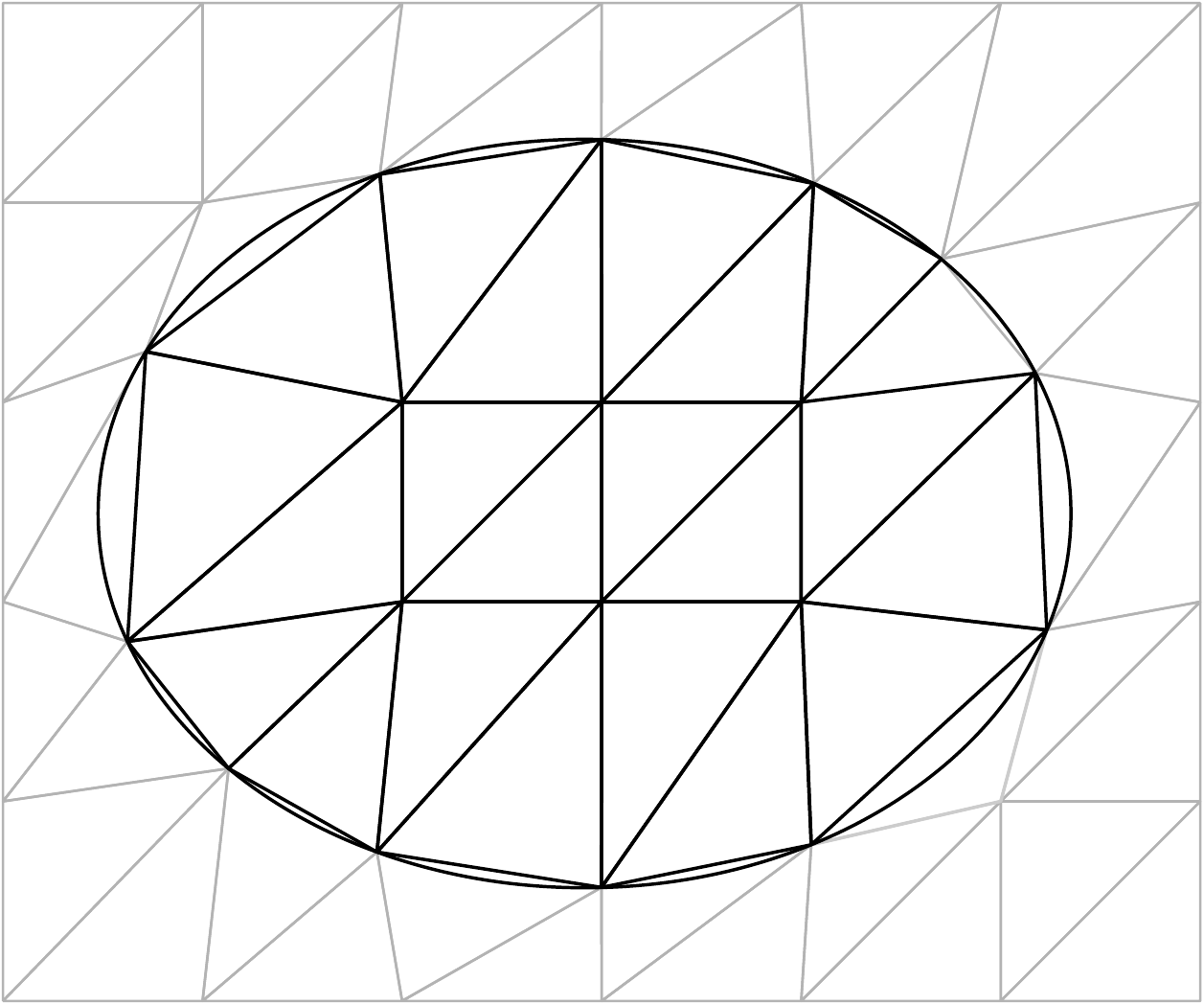}
		\caption{ }
		\label{fig:cfe-c}
	\end{subfigure}
	\caption{Grid generation: (a) the structured mesh $\widetilde{T}$ and the domain $\Omega$; (b) the adapted mesh $T^{\infty}$; (c) the active mesh $T_a$.}
	\label{fig:cfe}
\end{figure}

We assume the feasible area $\widetilde{\Omega}$ to be rectangular. $\widetilde{\Omega}$ is discretized by a regular triangular grid. The grid is denoted by $\widetilde{T}$, the number of elements by $\widetilde{N}^{el}$ and the number of nodes by $\widetilde{N}^{no}$. As in Figure \ref{fig:cfe-a}, the boundary $\delta\Omega_0$ of the shape to be optimized is superimposed onto the grid, represented, for example, by a set of splines.
In a second step, we adapt the grid $\widetilde{T}$ to the boundary of the shape by moving the closest vertex of a cut edge onto the boundary \cite{bolten2020}. The adapted grid is denoted by $T^{\infty}$. The computations are only performed on the elements inside the domain. We call these elements active elements, making up the active grid $T_{a}$, and the corresponding nodes are called active nodes ($P_{a}$, with $|P_{a}|=:N_a$). The active elements are a subset of the elements of $T^{\infty}$, which itself is a perturbation of $\widetilde{T}$. An important advantage of this approach is the possibility to implicitly build the system of equations for all nodes in $T^{\infty}$, with the rows and columns corresponding to the non-active nodes containing only zeros, and then perform the calculation only with the ``active" sub-mesh that we denote by $T_a$.

In this way, after updating the shape, the process starts again, each optimization step, with $\widetilde{T}$, with the same node numbering and connectivity. The mesh is updated according to \cite{bolten2020}. This gives a substantial speed-up in building the mesh in comparison with full re-meshing approaches, as only elements at the boundary have to be changed. Additionally, as most of the elements do not change, the system matrix has to be updated only for entries corresponding to nodes in elements that do change. Hence, the full assembly of the system matrix in each iteration is avoided. This provides an advantage over both the re-meshing and the mesh-morphing techniques, while obtaining an accuracy comparable to re-meshing approaches. 

\subsection{Mapping}
\label{sec:rminEG-map}
The structured meshing technique described in the previous section is decidedly well suited for designing a mapping of the type described in section \ref{sec:rminEG}.
 This is why, as a proof of concept, we introduce a quite simple mapping adjusted for this specific meshing technique, which still leads to a considerable speed up in many of our test cases.
Consider the linear systems $\MK(\rho^{(i)}) \Vu^{(i)} = \Vf^{(i)}$ and $\MK(\rho^{(i+1)}) \Vu^{(i+1)} = \Vf^{(i+1)}$ and $\MW^{(i)}\in\mathbb{R}^{N_a^{(i)}\times k},$ with $\textbf{range}(\MW^{(i)})$ approximating the invariant subspace corresponding to the $k$ smallest eigenvalues of $\MK(\rho^{(i)})$. As in this meshing approach the initial connectivity of the mesh is kept for all optimization steps, we can uniquely identify each node in iteration $(i+1)$ with a node in iteration $(i)$. We distinguish between three cases to determine the matrix $\widetilde{\MW}
^{(i)}$ representing the mapped approximate invariant subspace : (1) Matrix entries corresponding to inner nodes that stay inner nodes are kept; (2) matrix entries corresponding to nodes that change from inner node to boundary node or vice versa are recalculated according to (\ref{eq:map-matrix}); and (3) matrix entries corresponding to nodes that change from inactive to active are calculated in the following way. Assuming only small changes in the geometry, these nodes must be boundary nodes or close to boundary nodes. It is therefore likely that the values of former active nodes in their neighborhood provide a better approximation than a value resulting from interpolating former active and inactive nodes. We therefore propose for these points the following extrapolation:
Consider that the status of node $\tilde{x}_{\ell}^{(i+1)}$ changes from inactive to active. From the set of nodes that share a finite element, we consider only the subset of nodes that were active in iteration $(i)$. For each active node, $x_{\ell}^{(i+1)}$ we have such a set $S_a^{(i, \ell)}$. We set the entries of $\widetilde{\MW}^{(i)}$ corresponding to these nodes to the weighted mean of the values at $S_a^{(i, \ell)}$, given by
\begin{align}
(\widetilde{\MW}^{(i)})_{(\ell, j)}= \frac{1}{|S_a^{(i, \ell)}|-1}\sum_{s\in S_a^{(i, \ell)}} \frac{\sum_m||s_m-x_{\ell}^{(i+1)}||_2 -||s-x_{\ell}^{(i+1)}||_2}{\sum_m||s_m-x_{\ell}^{(i+1)}||_2} (\MW^{\infty,(i)})_{(s, j)},
\label{eq:extra}
\end{align}
with weights corresponding to the distance between the neighbors' location on the old grid and $x_{\ell}^{(i+1)}$.\\

\subsection{Test of Mapping an Approximate Invariant Subspace}
\begin{figure}
    \centering
    \includegraphics[width=0.8\textwidth]{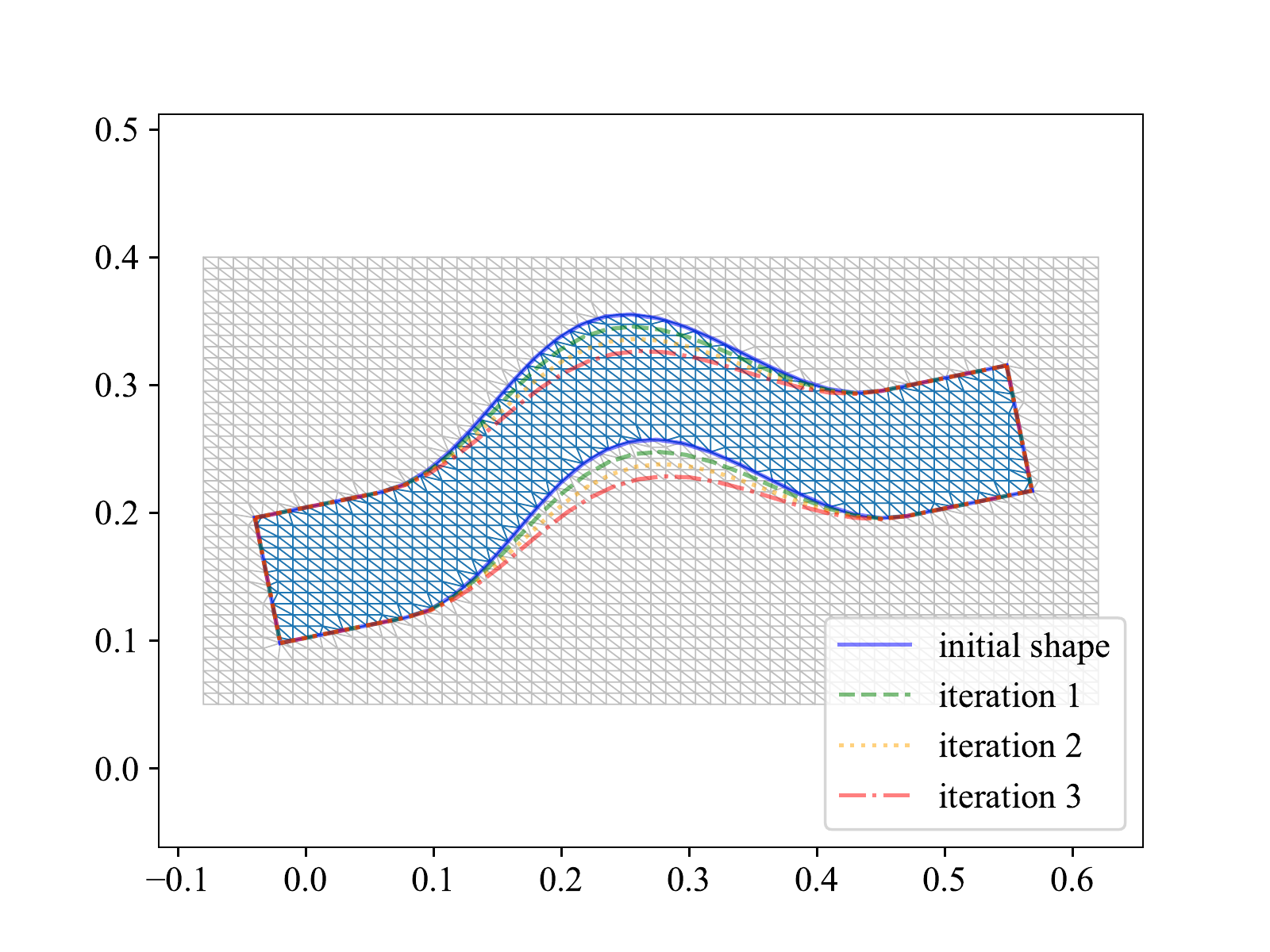}
    \caption{Sequence of shapes for a model problem.}
    \label{fig:bentrod}
\end{figure}

As our first example, we consider a model problem solving the linear elasticity equations on a bent rod with a $181\times121$ nodes grid, resulting in $5507$ active nodes. To get a first impression of the quality of the approximation of the invariant subspace via the mapping, we perform deformation steps of the shape in a controlled way, see Figure \ref{fig:bentrod}, and calculate the principal angles of the resulting mapped approximate invariant subspaces and the true invariant subspace corresponding to the smallest eigenvalues of the new system matrix, which has been computed for the purpose of comparison only.

\begin{figure}
    \centering
    \includegraphics[width=0.8\textwidth]{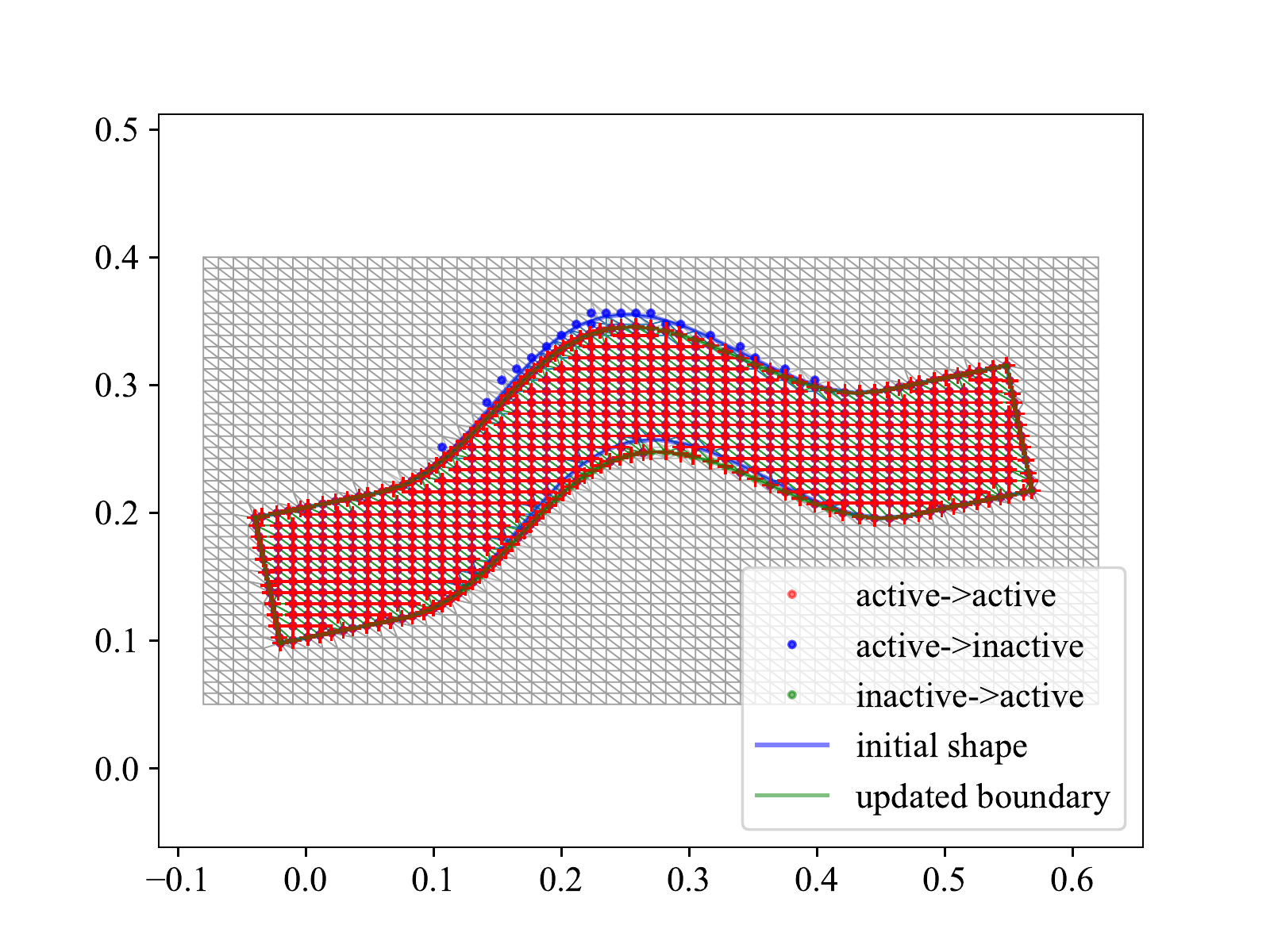}
    \caption{The status changes for nodes near the boundary for a small change in shape.}
    \label{fig:map-vis}
\end{figure}

We perform three steps of deformation of the original shape, with $N_{rec}=15$, and we calculate the principal angles between the mapped approximate invariant subspace and the invariant subspace corresponding to the 20 smallest eigenvalues according to \cite[p.~604]{golub1996}. We can see that most of the angles are rather small, i.e., our approach approximates these spaces quite well. Note that we only consider small geometric deformations.

  \begin{table}
  \small
		\begin{tabular}{l | *{10}{r}}
		no. & 1 & $\cdots$ & 8 & 9 & 10 & 11 & 12 & 13 & 14 & 15 \\
		\hline
		i 1 & 0.992 & $\cdots$ & 0.945 & 0.939 & 0.935 & 0.935 & 0.913 & 0.872 & 0.806 & 0.150\\
		i 2 & 0.983 & $\cdots$ & 0.904 & 0.843 & 0.785 & 0.722 & 0.551 & 0.331 & 0.153 & 0.102\\
		i 3 & 0.987 & $\cdots$ & 0.933 & 0.908 & 0.899 & 0.844 & 0.781 & 0.700 & 0.476 & 0.217
		\end{tabular}
		\caption{Cosines of principal angles (cos$(\theta_i )$) between the approximate invariant subspace and the true invariant subspace corresponding to the 20 smallest eigenvalues.}
		\label{tab:angles}
\end{table}

\section{Numerical Results}
\label{sec:num}

To demonstrate the efficacy of the described methods, we consider three two-dimensional examples. The first one concerns a turbine blade on which the Poisson equation is solved. For the second one, we solve the linear elasticity equation on a bent rod. Subsequently, we consider an example where Krylov subspace recycling does not work well using only the mapping described in section \ref{sec:rminEG-map}, and where it is necessary to apply the warm-start Krylov-Schur algorithm introduced in section \ref{sec:rminEG-kryschur}.

\subsection{Poisson equation on a turbine blade}
\label{sec:num-pois}
\begin{figure*}
	\centering
	\begin{subfigure}[t]{0.24\textwidth}
		\centering
		\includegraphics[width=\textwidth]{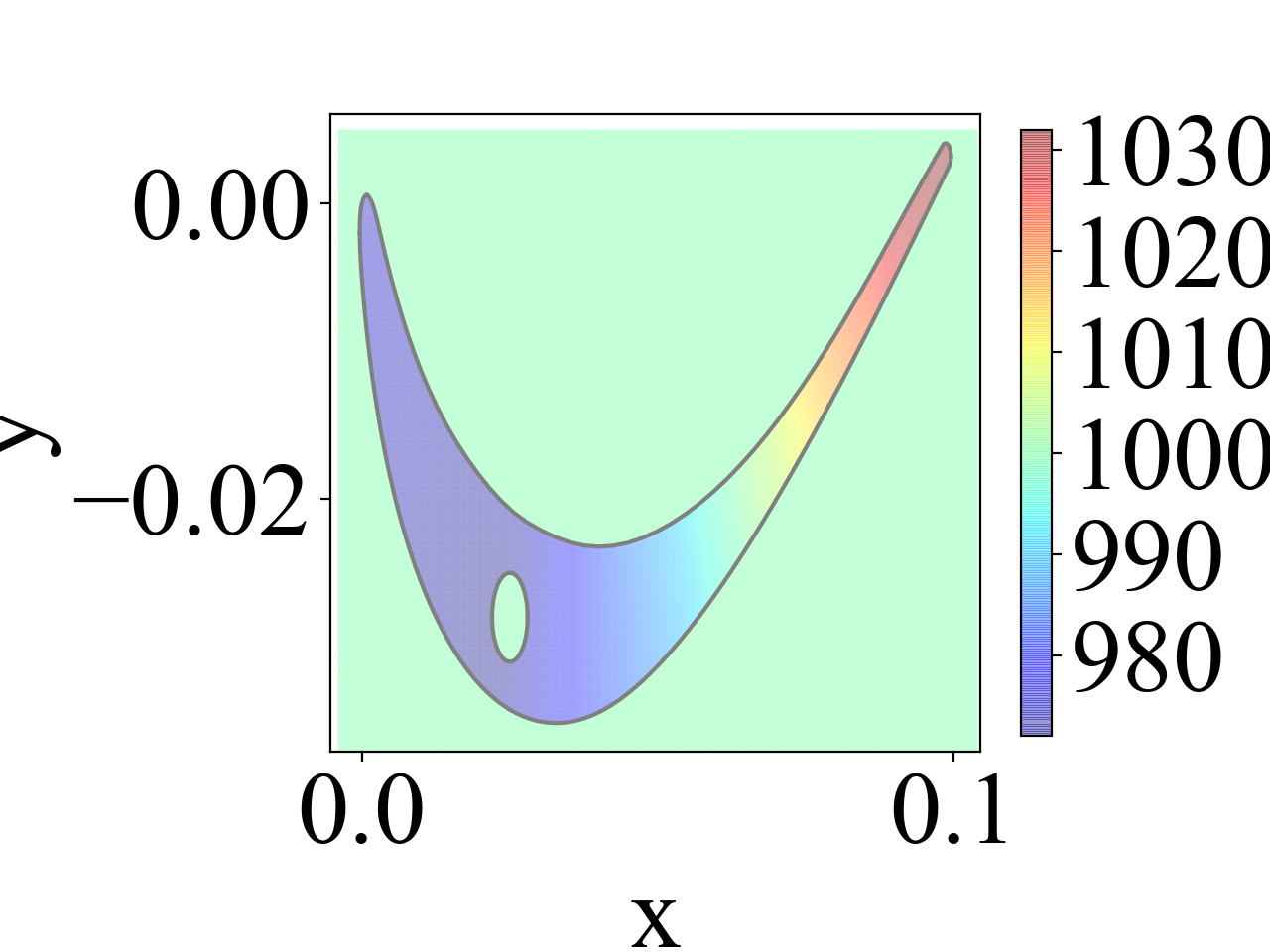}
		\caption[]%
		{{\small Initial configuration.}}    
		\label{fig:pois-blade1}
	\end{subfigure}
	\begin{subfigure}[t]{0.24\textwidth}  
		\centering 
		\includegraphics[width=\textwidth]{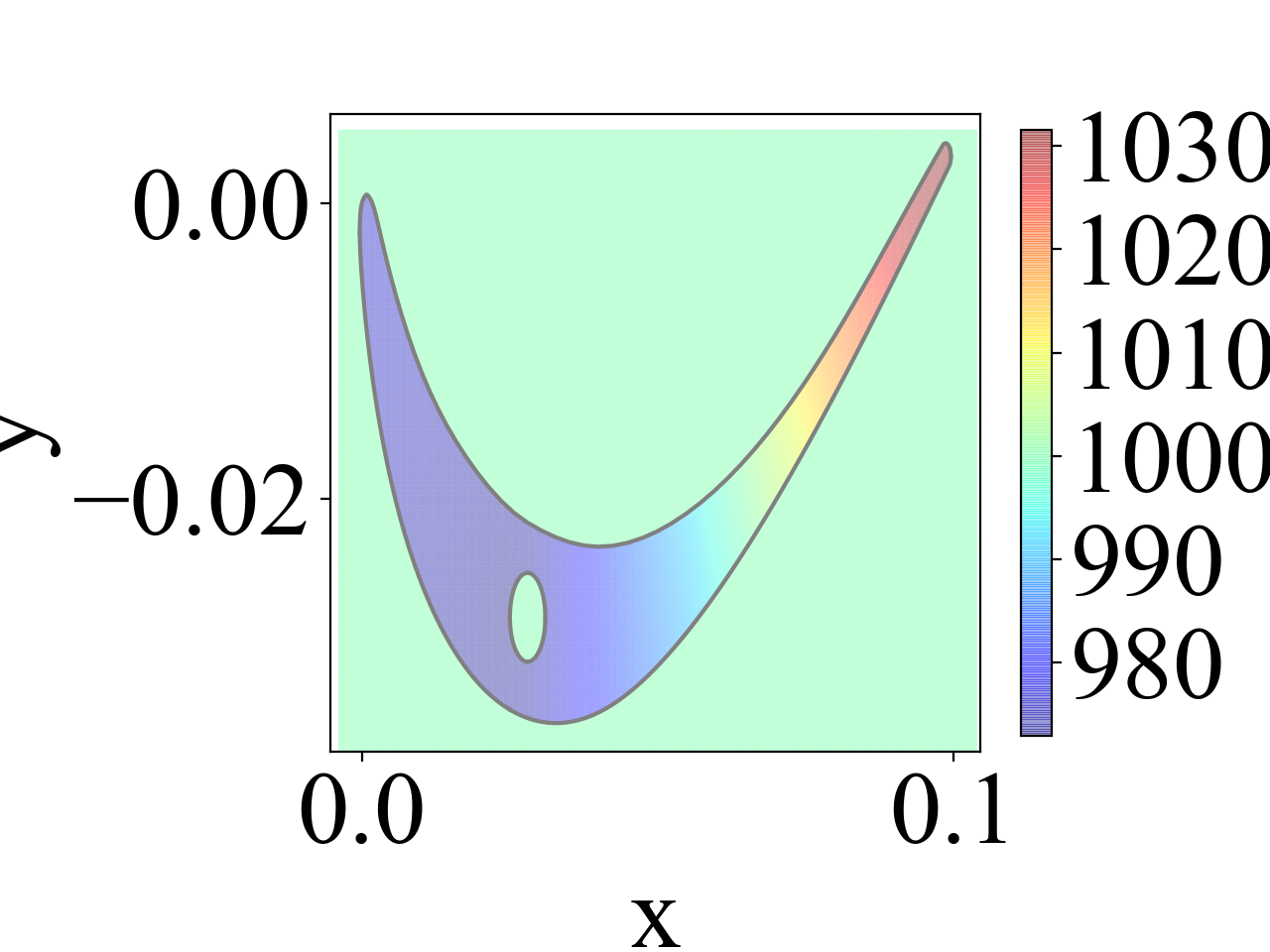}
		\caption[]%
		{{\small Configuration 1.}}    
		\label{fig:pois-blade2}
	\end{subfigure}
	\begin{subfigure}[t]{0.24\textwidth}   
		\centering 
		\includegraphics[width=\textwidth]{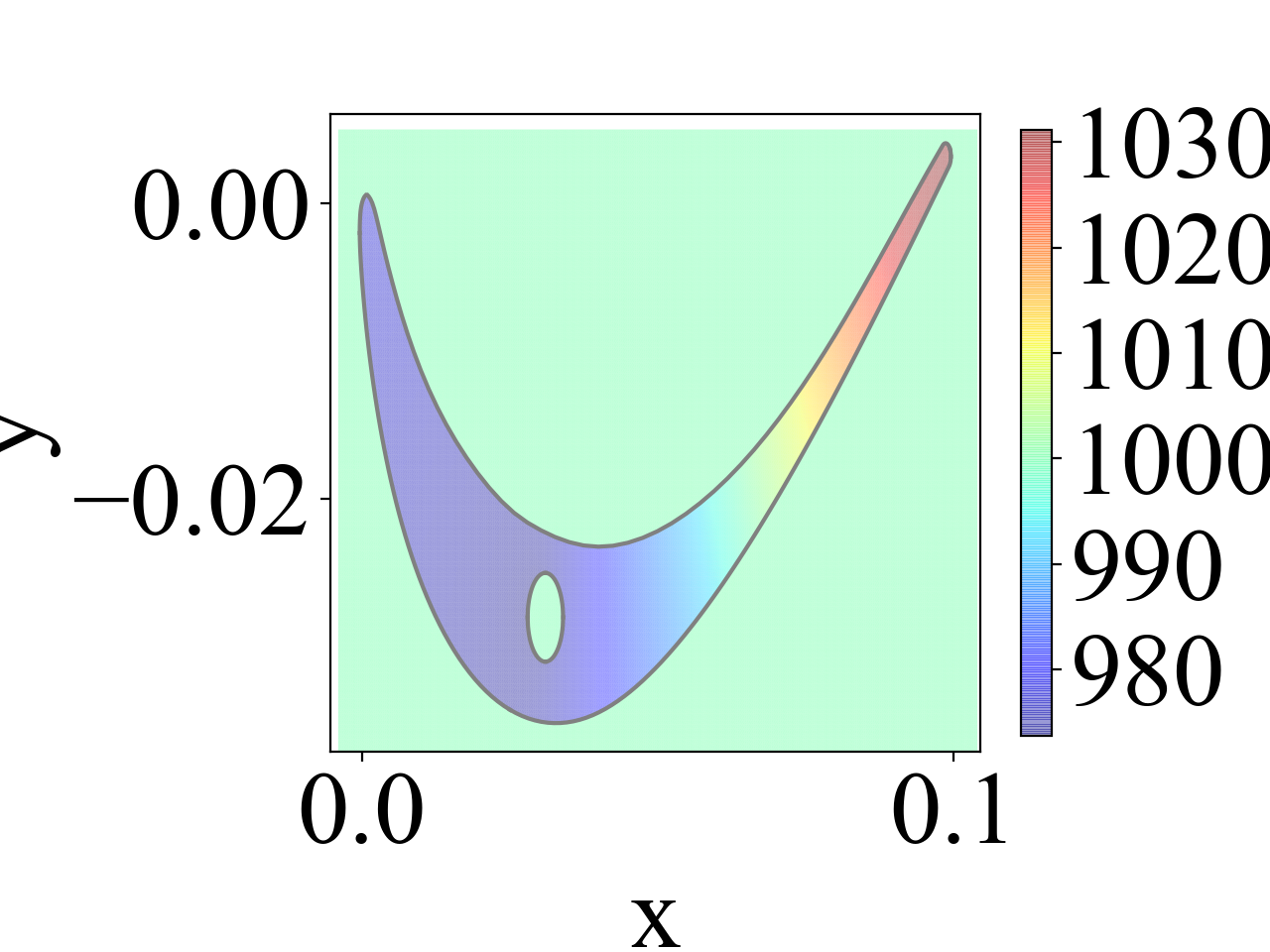}
		\caption[]%
		{{\small Configuration 2.}}    
		\label{fig:pois-blade3}
	\end{subfigure}
	\begin{subfigure}[t]{0.24\textwidth}   
		\centering 
		\includegraphics[width=\textwidth]{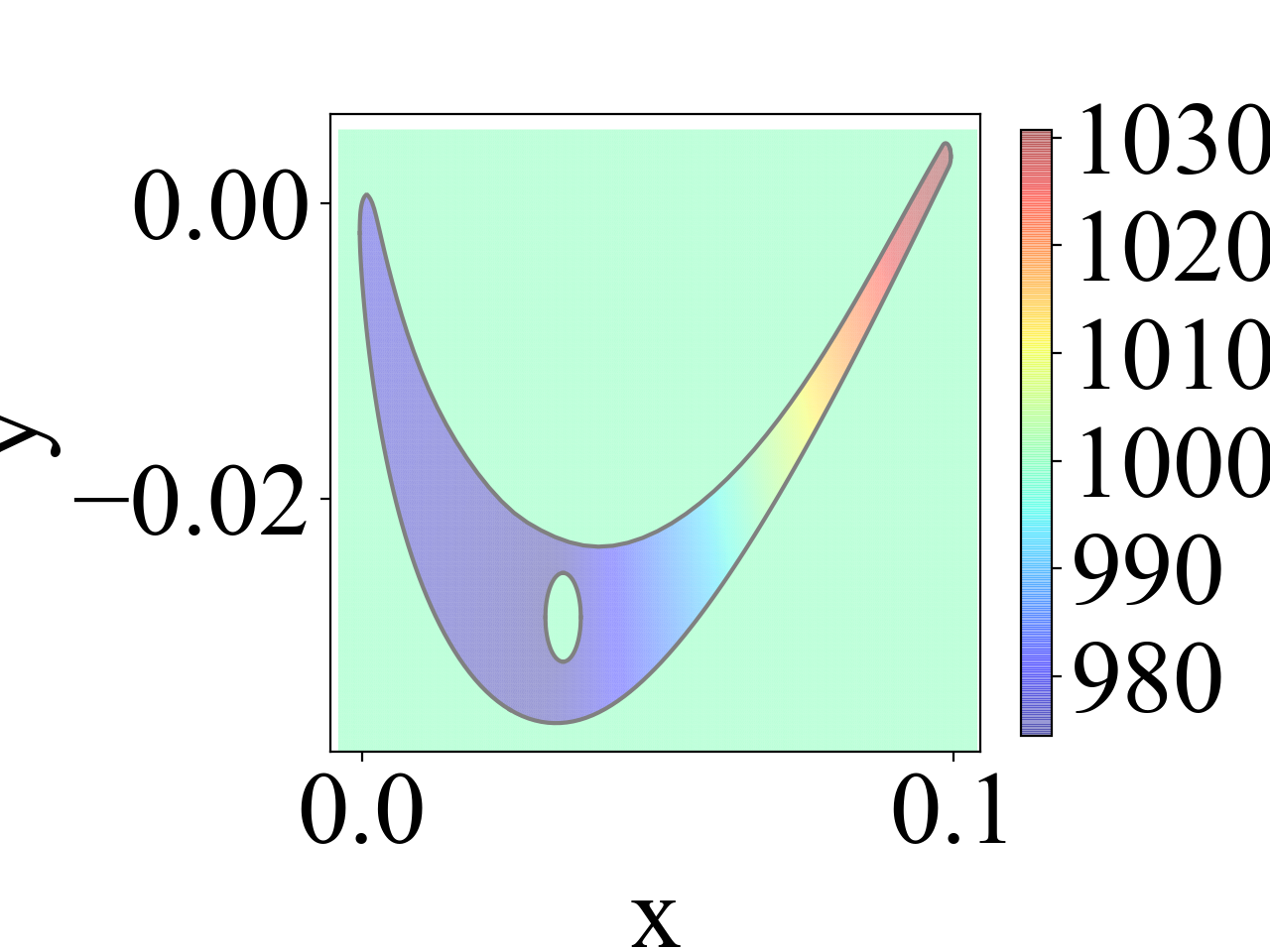}
		\caption[]%
		{{\small Configuration 3.}}    
		\label{fig:pois-blade4}
	\end{subfigure}
	\caption{\small Test problem 1: Poisson equation on a turbine blade. Three consecutive positions of the hole.} 
	\label{fig:pois-blade}
\end{figure*}

As a first example, we solve the Poisson equation on a turbine blade, where the change of geometry is caused by a hole representing a cooling channel as its position inside the turbine blade is being optimized. In this example, the movement of the hole is artificial and not driven by an optimization. $\widetilde{T}$ is a $361$ by $181$ grid, which gives close to $12,000$ active nodes on $T_a$. On the boundary of the blade, Robin-boundary conditions hold with constant heat coefficients, and the temperature on the cooling channel boundary is chosen to be two times lower than the one on the outer boundary. We perform an initial solve using MINRES and then three rMINRES solves for three consecutive changes of the domain; the geometries are visualized in Figure \ref{fig:pois-blade}. A simple $IC(0)$ preconditioner is used.

\begin{figure*}[htbp]
	\centering
	\begin{subfigure}[b]{0.53\textwidth}
		\includegraphics[width=\textwidth]{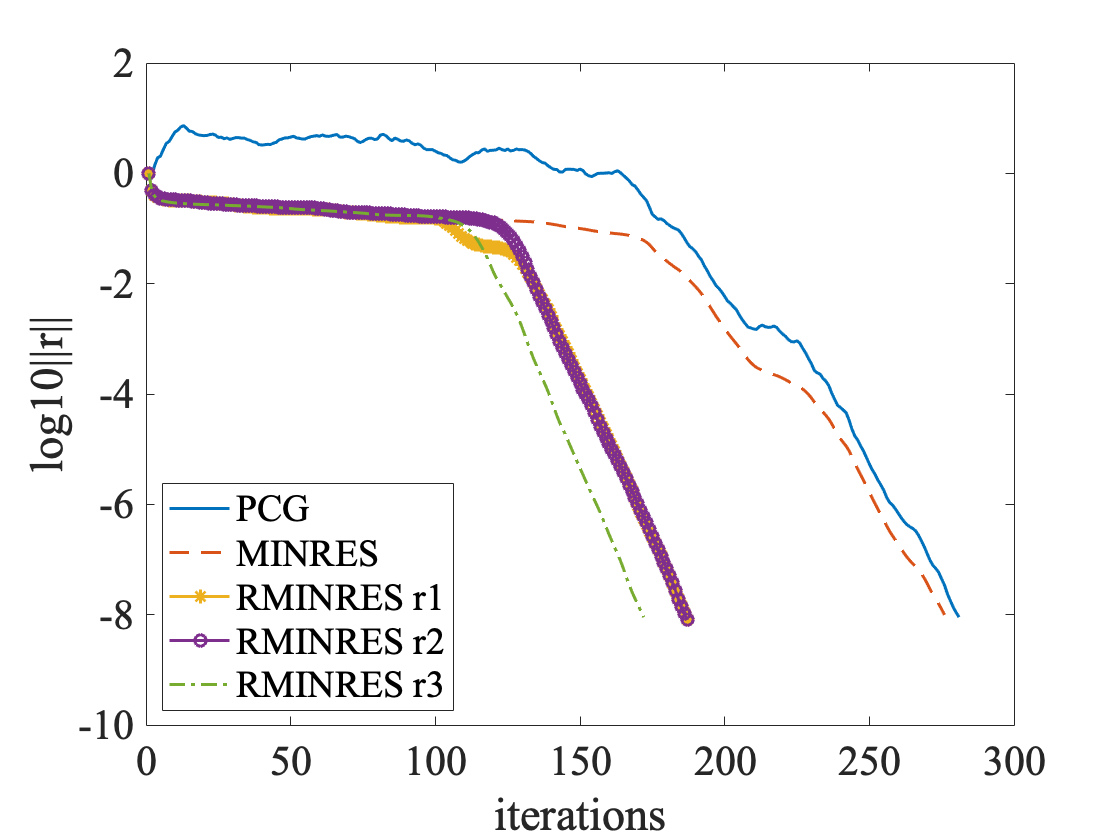}
		\caption{Residual norm convergence.}
	\end{subfigure}
	\begin{subfigure}[b]{0.45\textwidth}
	\begin{tabbing}
		\small
		\begin{tabular}{|p{5mm}||p{0.9cm}|p{0.8cm}|p{1.6cm}| }
			\hline
			\multicolumn{4}{|c|}{Data for $k=15$} \\
			\hline
			\small Opt. step & $N$ & \small 
			\# its & \small active \hspace{.5cm}  
			(inactive) \\
			\hline
			0   &  $11,893$   & 275 & -\\
			1&   $11,893$ & 186 & 210 (210)\\
			2 & $11,893$ & 186 & 210 (210)\\
			3    & $11,893$ & 171 & 210 (210)\\
			\hline
		\end{tabular}
	\end{tabbing}
\caption{Number of unknowns, $N$, number of rMINRES iterations, and number of nodes changing from active to inactive or vice versa.}
\end{subfigure}
\caption{Convergence results for solving Poisson's equation for 3 consecutive optimization steps.}
\label{fig:pois-data}
\end{figure*}

Figure \ref{fig:pois-data} demonstrates that a speed-up of more than 30\% is obtained. In this example, the number of nodes does not change, due to the fact that the hole is moved in very controlled way and without changes in size. Nevertheless, the rows and columns of the matrix do not represent the same nodes in the region of the hole.

\subsection{Gradient based shape optimization with linear elasticity as governing PDE}
\label{sec:num-le}

\begin{figure}
    \centering
    \includegraphics[width=0.7\textwidth]{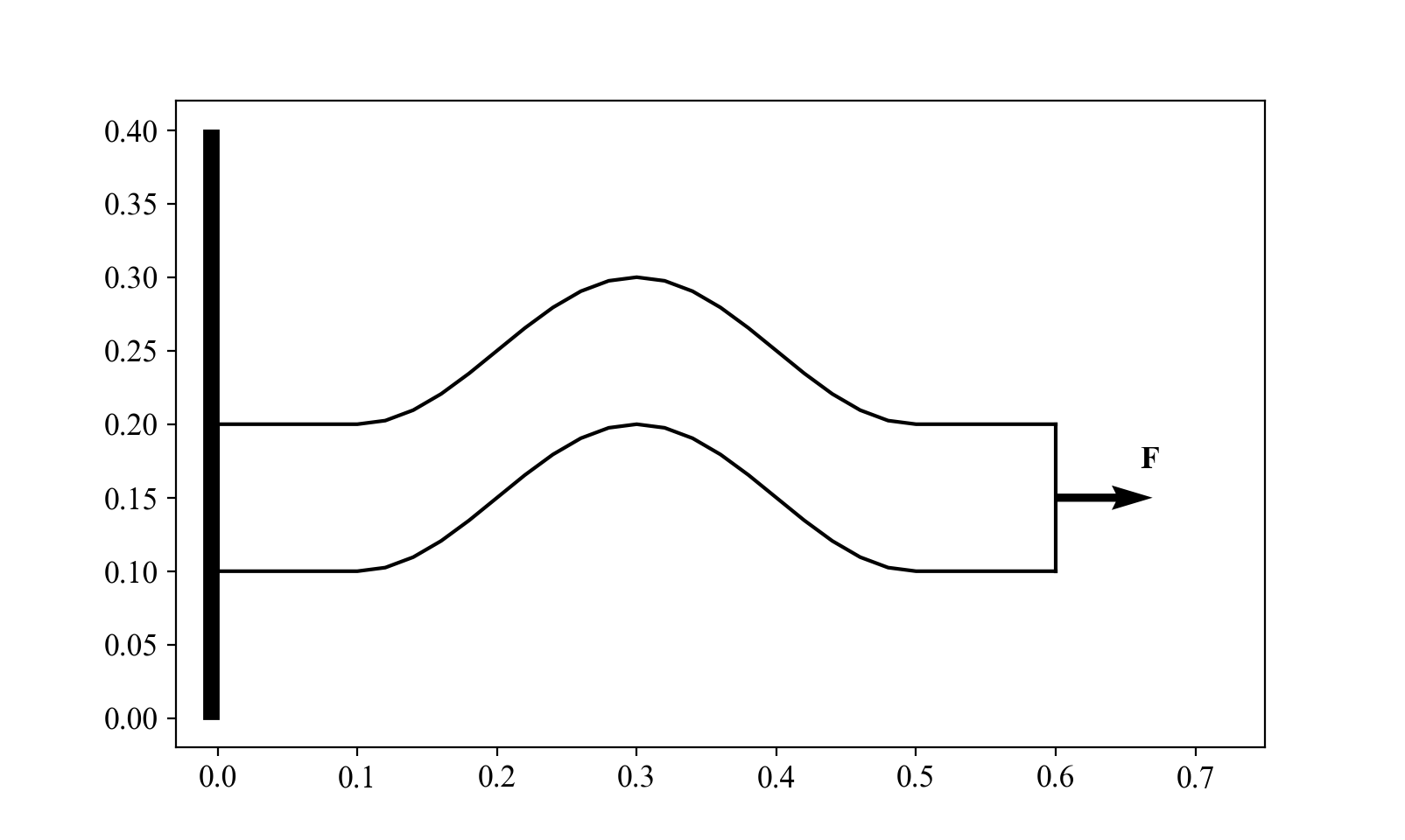}
    \caption{Test problem 2: Linear elasticity problem for ceramic object under tensile load.}
    \label{fig:LE-toy}
\end{figure}

\begin{figure*}
	\centering
	\begin{subfigure}[t]{0.24\textwidth}
		\centering
		\includegraphics[width=\textwidth]{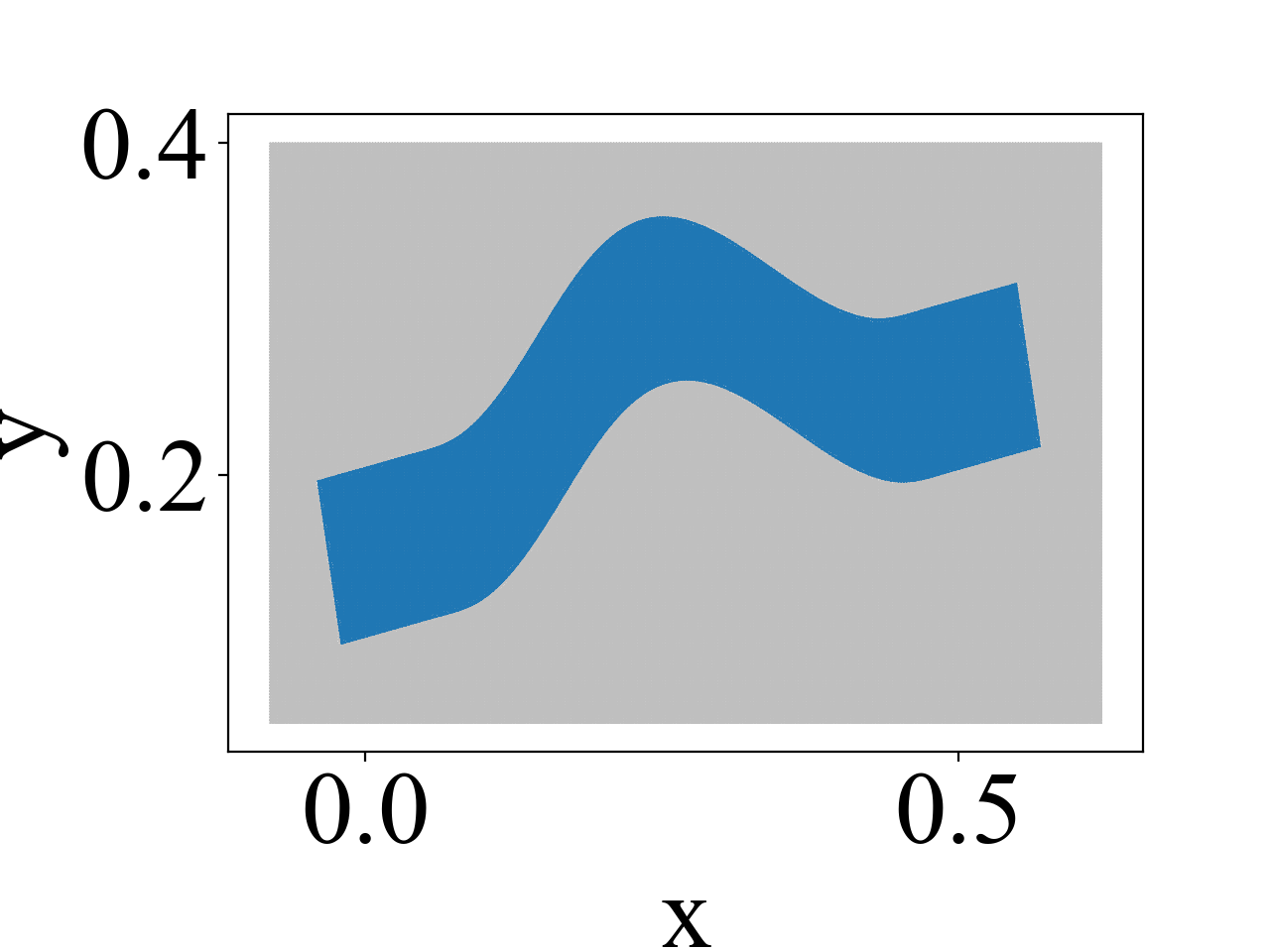}
		\caption[Network2]%
		{{\small Initial configuration.}}    
		\label{fig:LE-rod1}
	\end{subfigure}
	\begin{subfigure}[t]{0.24\textwidth}  
		\centering 
		\includegraphics[width=\textwidth]{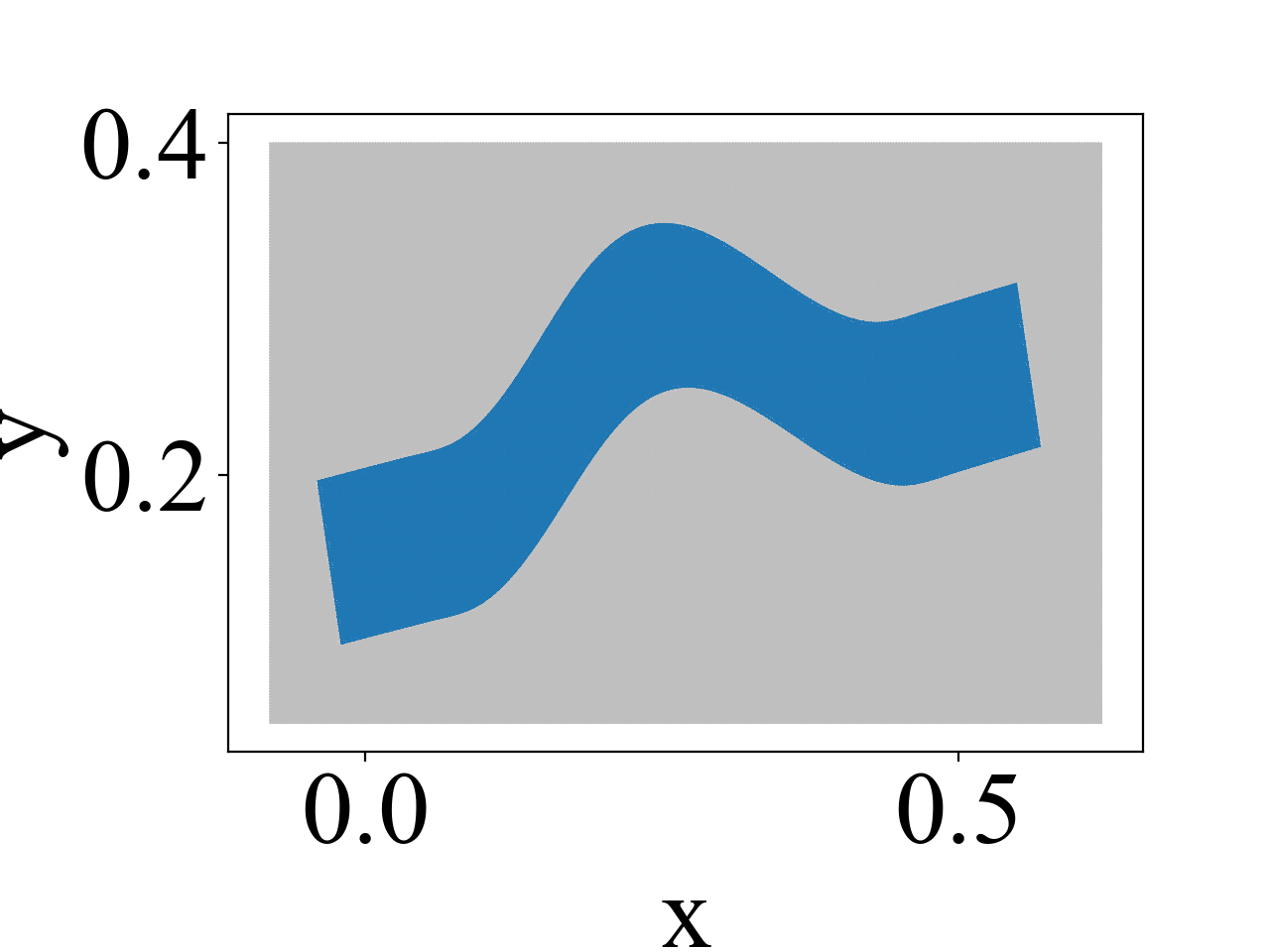}
		\caption[]%
		{{\small Configuration 1.}}    
		\label{fig:LE-rod2}
	\end{subfigure}
	\begin{subfigure}[t]{0.24\textwidth}   
		\centering 
		\includegraphics[width=\textwidth]{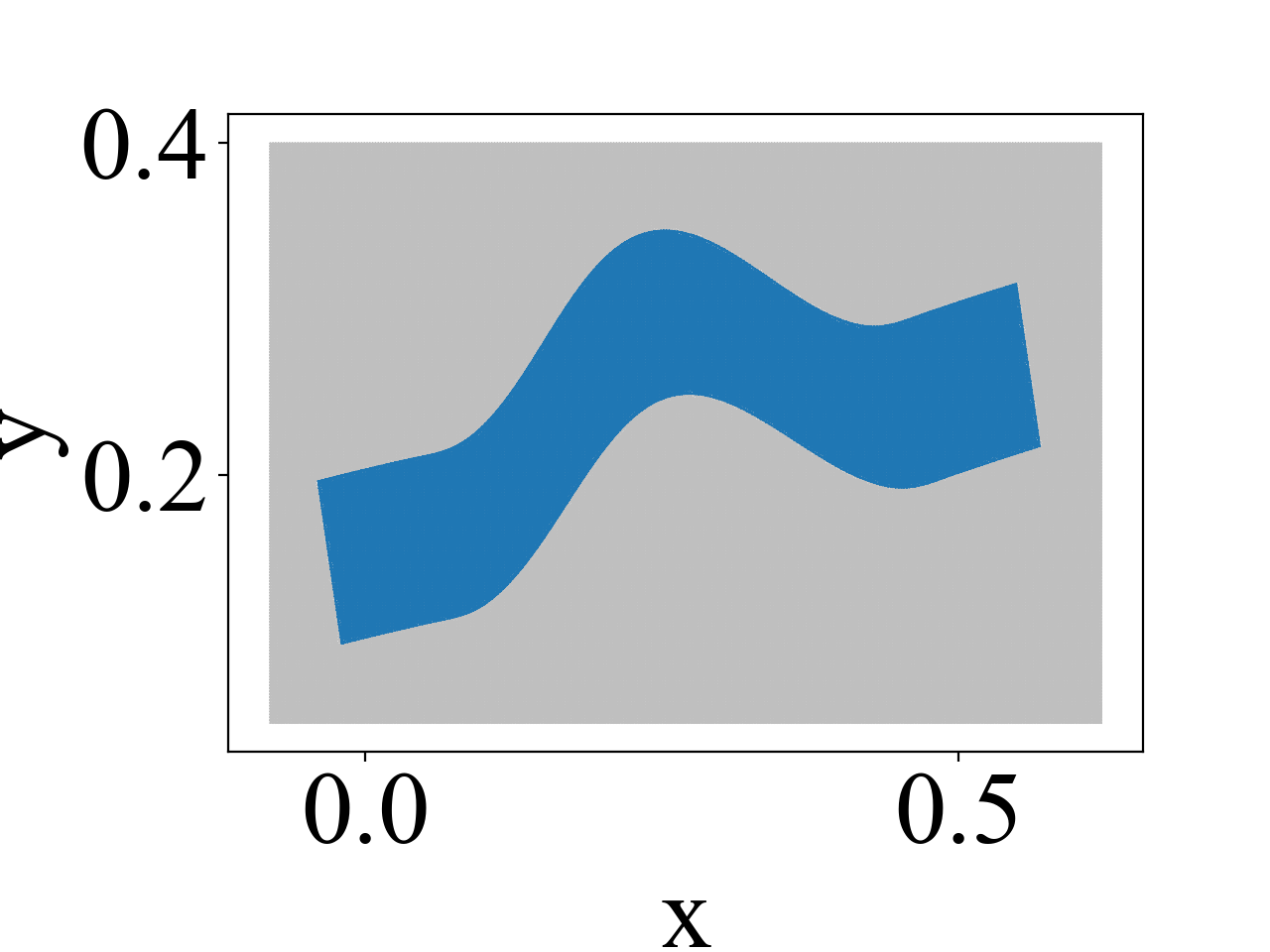}
		\caption[]%
		{{\small Configuration 2.}}    
		\label{fig:LE-rod3}
	\end{subfigure}
	\begin{subfigure}[t]{0.24\textwidth}   
		\centering 
		\includegraphics[width=\textwidth]{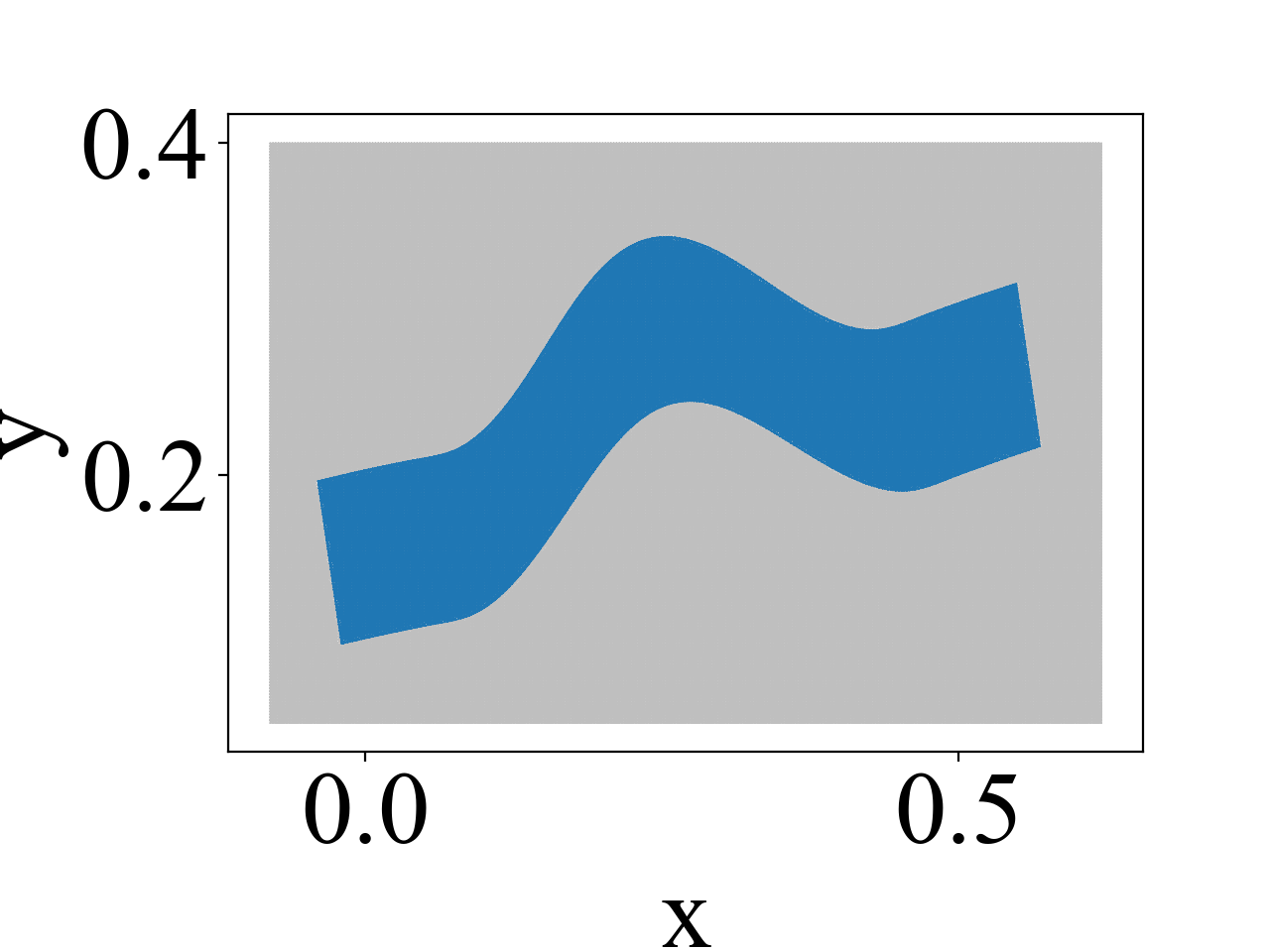}
		\caption[]%
		{{\small Configuration 3.}}    
		\label{fig:LE-rod4}
	\end{subfigure}
	\caption{\small Test problem 2: subsequent shapes in the optimization of a bent rod.} 
	\label{fig:LE-rod}
\end{figure*}

In our second example, we consider a bent rod that is clamped at the left, i.e., with zero-boundary conditions on the left boundary, and a tensile load is applied on the right boundary, i.e., Neumann-boundary conditions are applied on the right side; see Figure \ref{fig:LE-toy}. In this example, the deformations of the shape of the rod are not artificially generated but originate from an optimization procedure using shape derivatives \cite{haslinger2003}. We assume that the rod is made from a ceramic material, in this case $\text{Al}_2\text{O}_3$. Ceramic is a linear elastic material. The shape of the rod is being optimized to maximize its reliability under the applied tensile load. The reliability of the rod is measured by a functional giving the probability of failure of the rod under a given tensile load,
\begin{align}
\label{eqa:ObFun}
J(\Omega ,u):=\frac{\Gamma(\frac{d}{2})}{2\pi^{\frac{d}{2}}}\int\limits_{\Omega}\int\limits_{S^{d-1}}\left(\frac{(\Vn\cdot\boldsymbol{\sigma} (u)\Vn)^+}{\sigma_0}\right)^m d\Vn dx,
\end{align}
where $\Omega \subseteq \mathbb{R}^d$ is the domain,  $u\in H^1(\Omega,\mathbb{R}^d)$ is the solution of the governing linear elasticity equation, $B(u,v)=L(v)\text{, }\forall v\in H^{1}_{0}(\Omega,\mathbb{R}^d)$, $\sigma(u)$ is the stress tensor, $m\geq2$ is the Weibull modulus, and $\sigma_0$ is some positive constant. For simplicity, $m=2$ is assumed in this example.
The differentiability of the functional is shown in \cite{bolten2015}. We use the discrete adjoint method to calculate the shape derivative of the Lagrangian; for more details see \cite{Bolten2019}. $\widetilde{T}$ is a $301$ by $201$ grid that leads to approximately $15,000$ active nodes, and therefore about $30,000$ unknowns in the linear elasticity equation.
The calculation of the gradient is based on the Steklov-Poincaré type metric introduced in \cite{Schulz2016}. The rod is straightening during the iteration process, as visualized in Figure \ref{fig:LE-rod}. As the resulting linear systems are ill-conditioned, we allow for a drop tolerance of $tol=0.001$ for the IC-preconditioning. Additionally, the matrix is reordered via reverse Cuthill-McKee reordering \cite{George:1981}. As in the previous example, three rMINRES solves are performed after the initial MINRES solve on three consecutive configurations in the optimization process.

\begin{figure}[t!]
	\centering
	\begin{subfigure}[b]{0.53\textwidth}
		\includegraphics[width=\textwidth]{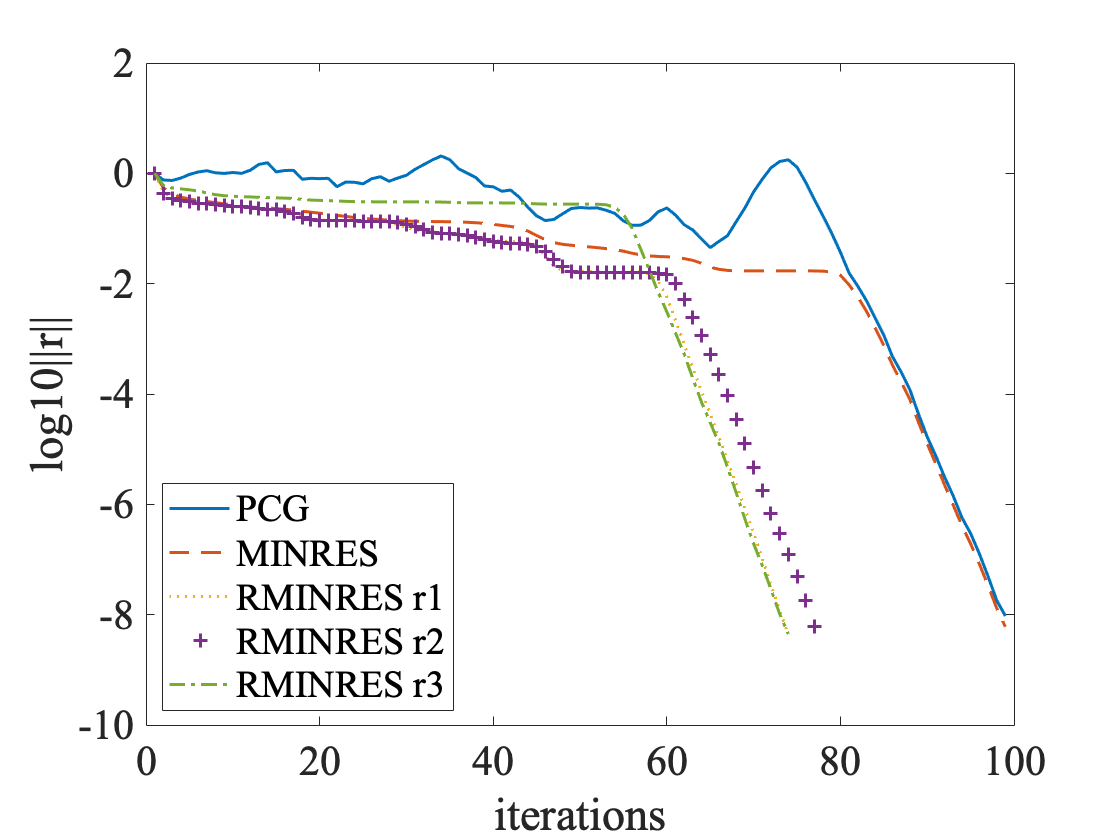}
		\caption{residual norm convergence}
	\end{subfigure}
	\begin{subfigure}[b]{0.45\textwidth}
		\begin{tabbing}
			\small
			\begin{tabular}{|p{5mm}||p{0.9cm}|p{0.8cm}|p{1.6cm}| }
				\hline
				\multicolumn{4}{|c|}{Data for $k=20$} \\
				\hline
				\small Opt. step & $N$ & \small \# its & \small active \hspace{0.3cm} 
				(inactive) \\
				\hline
	0   &  $30,108$   & $98$ & -\\
1&   $30,062$  & $73$ & $352$ $(375)$\\
2 & $30,030$ &  $76$ & $362$ $(387)$\\
3    & $29,982$ & $73$ & $353$ $(376)$\\
				\hline
			\end{tabular}
		\end{tabbing}
		\caption{Number of unknowns, $N$, number of rMINRES iterations, and number of nodes changing from active to inactive or vice versa.}
	\end{subfigure}
	\caption{Convergence results for solving the linear elasticity equation for 3 consecutive optimization steps.}
	\label{fig:le-data}
\end{figure}

The convergence is visualized in Figure \ref{fig:le-data}. Here we observe a speed up of the convergence of around $25\%$.

\subsection{Krylov-Schur enhanced recycling}
\label{sec:num-ks}

\begin{figure}[b!]
	\centering
	\begin{subfigure}[t]{0.32\textwidth}
		\centering
		\includegraphics[width=\textwidth]{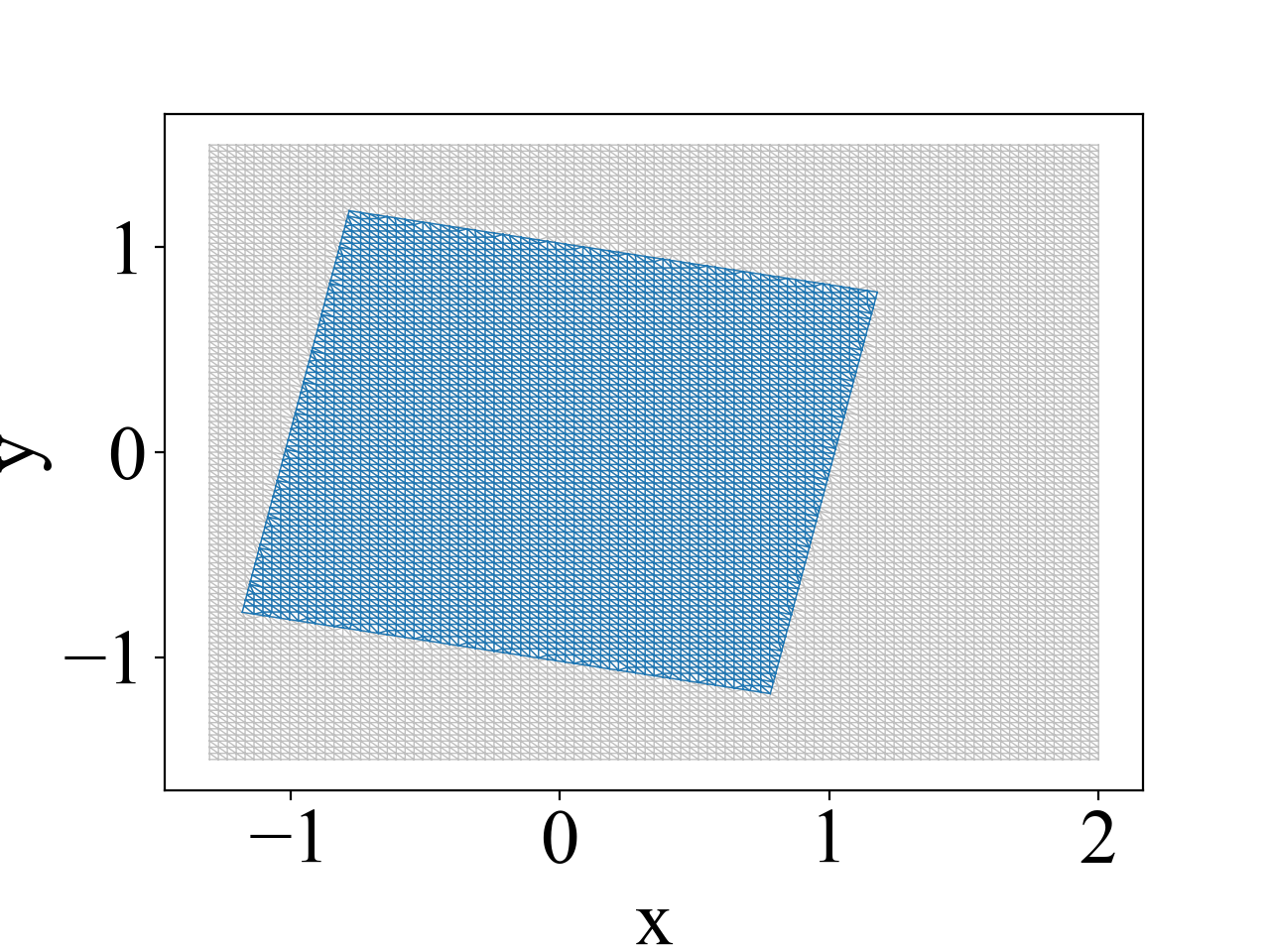}
		\caption[Network2]%
		{{\small Initial configuration.}}    
		\label{fig:ks1-geo1}
	\end{subfigure}
	\begin{subfigure}[t]{0.32\textwidth}  
		\centering 
		\includegraphics[width=\textwidth]{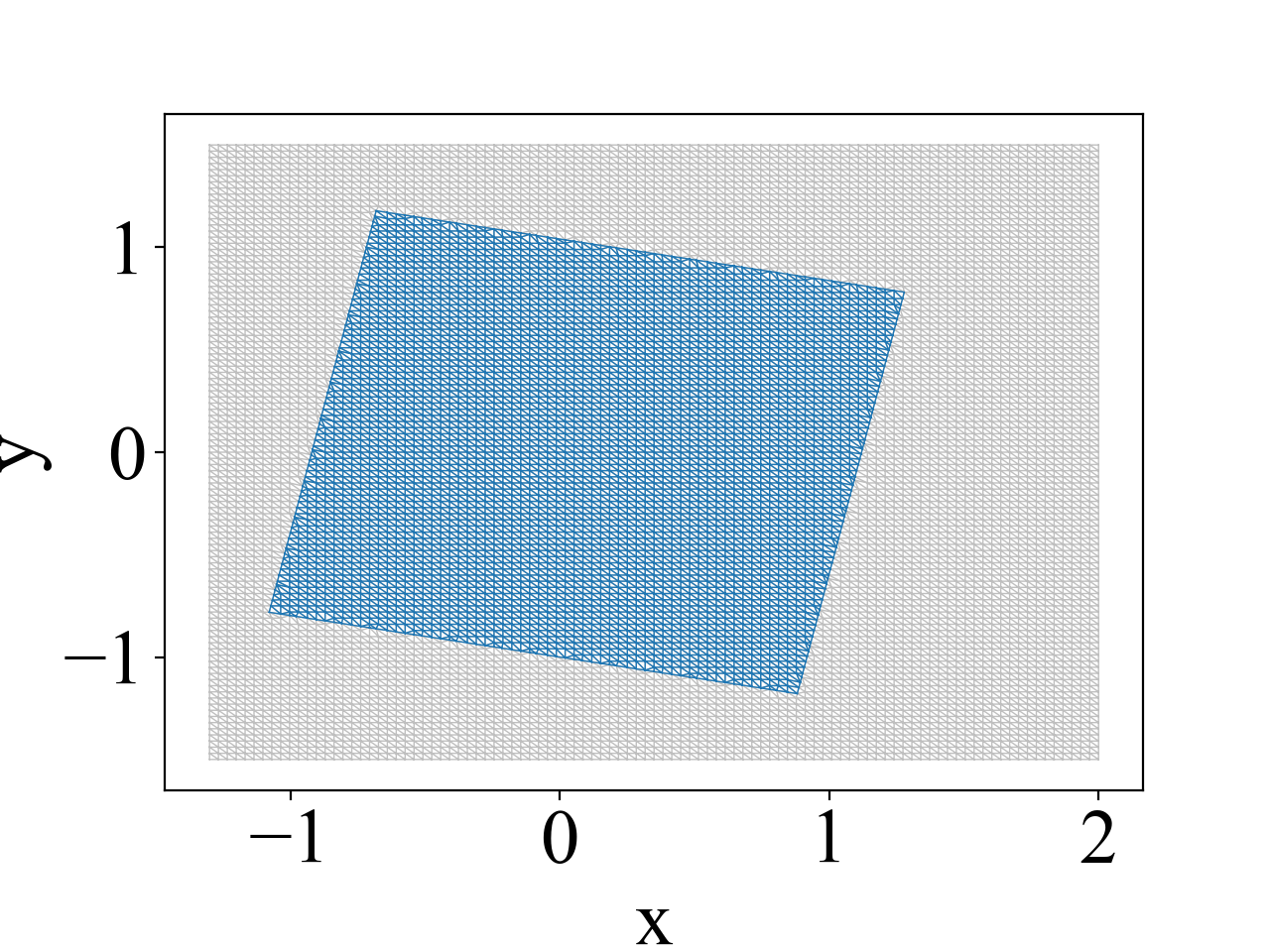}
		\caption[]%
		{{\small Configuration 1.}}    
		\label{fig:ks1-geo2}
	\end{subfigure}
	\begin{subfigure}[t]{0.32\textwidth}   
		\centering 
		\includegraphics[width=\textwidth]{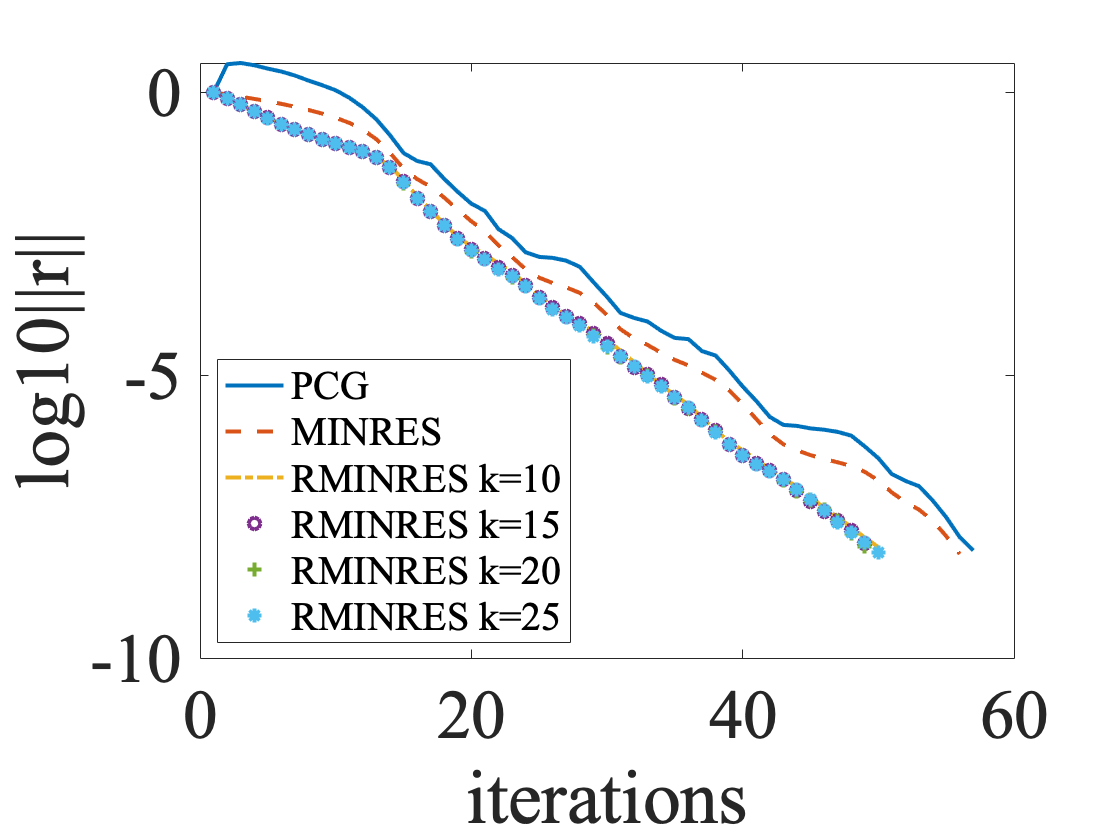}
		\caption[]%
		{{\small Residual norm convergence for several dimensions of the recycled space.}}    
		\label{fig:ks1-conv}
	\end{subfigure}
	\caption{\small Test problem 3: Example where  mesh-based mapping of the approximate invariant subspace is not sufficiently accurate.} 
	\label{fig:ks1}
\end{figure}

We have shown two examples where the proposed mapping of the approximate invariant subspace works quite well. However, in some cases the mapping
may fail to produce a sufficiently accurate approximate invariant subspace.

We consider the classical example for the Poisson equation on a square with source term $f(x)=1$ and zero-boundary conditions. The square is rotated and placed somewhere inside the feasible region. We discretize the feasible area by a $101\times 101$ grid. The square is not deformed but moved along the $x$-axis by a distance of $0.1$, as visualized in Figure \ref{fig:ks1}.

\begin{figure}[t!]
    \centering
    \includegraphics[width=0.7\linewidth]{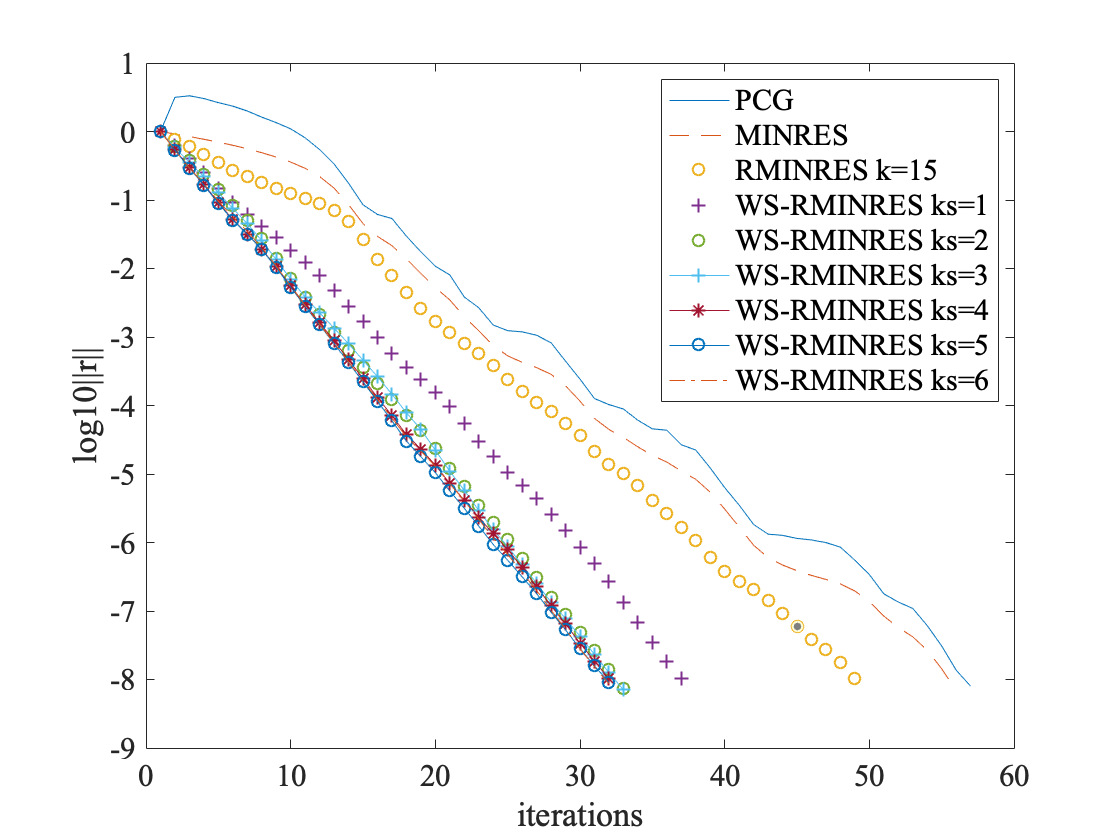}
    \caption{Convergence results with Krylov-Schur enhanced approximate invariant subspaces for $k=15$.}
    \label{fig:ks2}
\end{figure}

Although the shape hasn't actually changed, in a Hausdorff sense the two domains and thus the corresponding grids representing the geometries are too far away from each other, so that we find a situation where an additional improvement of the approximate invariant subspace,
as described in section \ref{sec:rminEG-kryschur},
is required for fast rMINRES convergence.
In Figure \ref{fig:ks1-conv} we show that the recycle space used  fails to lead to a considerable speed up, even if we increase the dimension of the recycle space $k$ from $10$ to $15$, $20$, or even $25$. We therefore apply the proposed warm-start Krylov-Schur algorithm with $tol=2 \cdot 10^{-8}$ and we vary the number of restarts of the warm-start Krylov-Schur algorithm from $1$ to $6$ to analyze the effectiveness.

Figure \ref{fig:ks2} shows the convergence  of rMINRES for each number of restarts of the warm-start Krylov-Schur algorithm performed before executing the rMINRES. It shows that, in this example, already two restarts of the Krylov-Schur reduces the iteration count from $48$ to $32$.

\begin{table}[t]
    \centering
    \small
    \begin{tabular}{c | *{3}{c}}
		angle & $\widetilde{W}$ & $U1$ & $U2$\\
		\hline
       1 & 0.99999844 & 0.99999837 & 0.99999849 \\
       2 & 0.99987744 & 0.99997532 & 0.99998698 \\
        3 & 0.99926389 &  0.99986194 & 0.99975142 \\
       4 & 0.99911852 &  0.99982407 & 0.99971646 \\
       5 & 0.99744298 &  0.99939602 & 0.99958597 \\
        6 & 0.99601106 &  0.9992919 & 0.99929083 \\
        7 & 0.98993238 & 0.99915824 & 0.99900528 \\
        8 & 0.89092459 &  0.99764084 & 0.9987683 \\
        9 & 0.86712723 &  0.99711616 & 0.99846744 \\
        10 & 0.72121063 &  0.99516303 & 0.99790027 \\
        11 & 0.29101069 & 0.99440875 & 0.9969349 \\
        12 & 0.12916091 &  0.99240021 & 0.99475094 \\
        13 & 0.06441871 & 0.98388291 & 0.99099253 \\
        14 & 0.02498528 & 0.93077594 & 0.99025971 \\
        15 & 0.00743539 &  0.71541515 & 0.95873314 \\
        16 &  -         &  0.14455214 & 0.29849549
    \end{tabular}
    \caption{Cosines of principal angles (cos$(\theta_i )$) between the invariant subspace corresponding to the smallest 20 eigenvalues and the approximate invariant subspaces from the mesh based mapping, $\widetilde{W}$, and after $1$
    restart, $U1$, 
    and $2$ restarts, $U2$, of the warm-start Krylov-Schur algorithm.}
    \label{tab:angles-ks}
\end{table}

To analyze the effect of the Krylov-Schur algorithm, we look at the principal angles between the invariant subspace corresponding to the $20$ smallest eigenvalues of $\MK$ and the respective approximate invariant subspaces in Table \ref{tab:angles-ks}. We can see that the approximation quality of the subspace improves significantly after one restart already.

\section{Conclusion}
\label{sec:conc}
In this paper, we have introduced new approaches to recycle information from the Krylov subspaces of previous systems in shape optimization. In contrast to other approaches previously considered in the literature, in shape optimization we have to face a changing number of unknowns, and the mapping from the unknowns in the algebraic systems to the meshes may not be consistent from one optimization step to the next. This makes it difficult to map subspaces from one optimization step to the next. In addition, the subspaces may not be sufficiently accurate, e.g., when the domain is moved. The change in the number of unknowns can be dealt with by evaluating a function represented by the coefficients of the solution of the previous optimization step at the current node location, possibly in combination with extrapolation if a new active node point was outside of the previous active mesh. If the approximate invariant subspace obtained from the previous iteration step is not accurate enough after having been transferred to the current geometry, a proposed warm-start Krylov-Schur algorithm can be used to improve the accuracy of this subspace.

The numerical results demonstrate the efficacy of our finite element-based approach in two different cases: (1) When the number of unknowns does not change due to small changes in the domain, only, and (2) when the number of unknowns does change due to larger deformations of the domain under consideration. The usefulness of the warm-start Krylov-Schur-based approach has been confirmed if the mesh-based mapping does not produce a sufficiently accurate approximate invariant subspace. In all of these cases the number of iterations necessary to solve the system has been substantially reduced.

\section*{Funding}
This work was supported by the federal ministry of research and education of Germany (BMBF, grant-no: 05M18PXA) as a part of the GIVEN consortium. This material is based  upon work supported by the National Science Foundation under Grant No. 1720305. 

\bibliography{literature}
\end{document}